\documentclass{article}

\usepackage{amsmath}
\usepackage{amsthm}
\usepackage{amsfonts}
\usepackage{amssymb}
\usepackage{upref}
\usepackage{srcltx}

\theoremstyle{plain}
\newtheorem{theorem}{Theorem}[section]
\newtheorem{proposition}[theorem]{Proposition}
\newtheorem{lemma}[theorem]{Lemma}

\theoremstyle{remark}
\newtheorem{remark}[theorem]{Remark}

\theoremstyle{definition}
\newtheorem{definition}[theorem]{Definition}

\numberwithin{equation}{section}

\newenvironment{mlist}{\list{}{\listparindent 0pt
\itemsep 0pt \parsep  0pt \topsep 0pt
}}{\endlist}

\newcommand{\eqdef}{\stackrel{\mathrm{def}}{=}}
\newcommand{\ad}{\operatorname{ad}}
\newcommand{\Ad}{\operatorname{Ad}}
\newcommand{\ddim}{\operatorname{ddim}}
\newcommand{\grad}{\operatorname{grad}}
\newcommand{\rank}{\operatorname{rank}}

\newcommand{\bbC}{\mathbb{C}}
\newcommand{\bbR}{\mathbb{R}}

\newcommand{\fg}{\mathfrak{g}}

\newcommand{\fk}{\mathfrak{k}}
\newcommand{\fm}{\mathfrak{m}}
\newcommand{\fz}{\mathfrak{z}}

\newcommand{\ft}{\mathfrak{t}}

\begin{document}

\title{\bf Integrability
of geodesic flows for metrics on suborbits of the adjoint
orbits of compact groups}
\author{Ihor V. Mykytyuk \\ \normalsize
email: mykytyuk{\_\,}i@yahoo.com
}

\maketitle

\begin{abstract}\noindent
Let $G/K$ be an orbit of the adjoint representation of a compact
connected Lie group
$G$, $\sigma$ be an involutive automorphism of $G$ and
$\tilde G$ be the Lie group of fixed points of
$\sigma$. We find a sufficient condition for the complete
integrability of the geodesic flow of the Riemannian metric on
$\tilde G/(\tilde G\cap K)$ which is induced by the
bi-invariant Riemannian metric on
$\tilde G$. The integrals constructed here are real analytic
functions, polynomial in momenta. It is checked that this
sufficient condition holds when
$G$ is the unitary group $U(n)$ and $\sigma$
is its automorphism determined by the complex conjugation.
\end{abstract}

\section*{Introduction}
\footnote{This research is partially supported by
Ministry of Economy and Competitiveness, Spain - CSIC,
under Project MTM2011-22528.}

Let $G/K$ be a homogeneous space of a compact Lie group
$G$. We consider the problem of the complete integrability
of the geodesic flow of the Riemannian metric on
$G/K$ which is induced by a bi-invariant Riemannian metric
on $G$. This problem was solved for some types of homogeneous
manifolds including symmetric spaces, spherical spaces,
Stiefel manifolds, flag manifolds, orbits of the adjoint
actions and others (see~~\cite{Ma78}, \cite{Mi82},
\cite{GS83}, ~\cite{My86}, \cite{MS00}, \cite{BJ01},
\cite{BJ04}, \cite{MP04}). Here we consider the new family
of homogeneous manifolds, namely, the suborbits of orbits of
the adjoint actions.

This paper is motivated by the paper~\cite{DGJ09}, in which
there were constructed integrable geodesic flows of
$\tilde G$-invariant metrics on the homogeneous space
$\tilde G/\tilde K=SO(n)/\bigl(SO(n_1)\times\cdots\times
SO(n_p)\bigr)$,
$n_1+\cdots+n_p=n$ ($n_i\geqslant1$,
$p\geqslant3$). The method of the proof in~\cite{DGJ09} is
based on investigations of bi-Poisson structures on the Lie algebras
$\mathfrak{u}(n)$ and $\mathfrak{so}(n)$
associated with Lie algebra deformations. We consider the
Lie-algebraic aspects of the integrability problem for such
homogeneous spaces. Our approach is based on the following
observation: the space
$\tilde G/\tilde K$ is a
$\tilde G$-suborbit of the adjoint orbit
$G/K=U(n)/\bigl(U(n_1)\times\cdots\times U(n_p)\bigr)$
of the Lie algebra $\mathfrak{u}(n)$ of the unitary group, i.e.
$\tilde G/\tilde K=\Ad(\tilde G)(a)$, where
$a\in\mathfrak{u}(n)$, and
$G/K=\Ad(G)(a)$. Moreover, $\tilde G$
is the group of fixed points of the involutive automorphism
$\sigma$ of $U(n)$ induced by the complex conjugation. In
other words, the space
$\tilde G/\tilde K$ is determined uniquely by the pair
$(G/K,\sigma)$, where $G/K=\Ad(G)(a)$
is an arbitrary adjoint orbit of the Lie group
$G$ in its Lie algebra $\fg$ with $a\in(1-\sigma_{*})\fg$,
$\sigma_*$ is the tangent automorphism of the Lie algebra
$\fg$.

Let $G$ be an arbitrary compact connected Lie group
$G$ with an involutive automorphism $\sigma:G\to G$ and let
$\tilde G$ be the set of fixed points of
$\sigma$. In the article we investigate the integrability of
the geodesic flow on the cotangent bundle
$T^*(\tilde G/\tilde K)$ defined by a
$\tilde G$-invariant metric on $\tilde G/\tilde K$
which is induced by a bi-invariant Riemannian metric on
$\tilde G$. As a homogeneous space
$\tilde G/\tilde K$, we consider the homogeneous space
associated with the adjoint orbit
$G/K=\Ad(G)(a)$ of an arbitrary point
$a\in(1-\sigma_{*})\fg$, i.e.
$\tilde K=\tilde G\cap K$. We find a sufficient purely
algebraic condition for the integrability of this geodesic
flow on the symplectic manifold
$T^*(\tilde G/\tilde K)$ (Theorem~\ref{th.10},
Propositions~\ref{pr.14} and~\ref{pr.15}). We prove that
this sufficient condition holds when
$G$ is the unitary group $U(n)$ and $\sigma$
is its automorphism defined by the complex conjugation
(Theorem~\ref{th.16}). These results confirm the
Mishchenko-Fomenko conjecture~\cite{MF78a}, which says that
a noncommutative integrability implies a Liouville
(complete) integrability by means of the integrals in the
same functional class (the geodesic flow of a metric on
$G/K$ induced by a bi-invariant metric on $G$
is integrable in a noncommutative sense by means of
analytic, polynomial in momenta functions~\cite{BJ01}). Our
approach is based on the fact that
$\tilde G/\tilde K\subset G/K$ is a totally real
(Lagrangian) submanifold of the homogeneous K\"ahler
manifold (the compact orbit)
$G/K$ and $T(\tilde G/\tilde K)\subset T(G/K)$
is a totally real submanifold of
$T(G/K)$. However, to simplify calculations we reformulate
this fact in some algebraic terms (not explicitly because
explicit reformulation would be very complicated from the point of
view of calculations on $T(T(G/K))$).
Also matters are considerably simplified by using
a very interesting observation of Bolsinov and Jovanovi\'c~\cite{BJ04}
about the set of roots associated with the homogeneous
space (orbit) $G/K$ and results of Kostant~\cite{Ko63}
describing all maximal dimension orbits of
semisimple complex Lie algebras.

One calls a Hamiltonian system on
$T^*M$ (completely) integrable if it admits the maximal
number of independent integrals in involution, i.e.
$\dim M$ functions commuting with respect to the Poisson
bracket on
$T^*M$ whose differentials are independent in an open dense
subset of $T^*M$. This set of integrals is a complete involutive
subset of the algebra
$C^\infty(T^*M)$. By Liouville's theorem the integral curves
of an integrable Hamiltonian system under a certain
additional compactness assumption are quasiperiodic (are
orbits of a constant vector field on an invariant torus).

Let $A^G$ be the set of all $G$-invariant
real analytic functions on the cotangent bundle
$T^*M$ of $M=G/K$. This space is an algebra with respect to
the canonical Poisson bracket on the symplectic manifold
$T^*M$. The natural extension of the action of
$G$ on $M$ to an action on the symplectic manifold
$T^*M$ is Hamiltonian with the moment mapping
$\mu^\mathrm{can}: T^*M\to\fg^*$. The functions of the type
$h\circ \mu^\mathrm{can}, h:\fg^*\to\bbR$,
are integrals for any
$G$-invariant Hamiltonian flow on $T^*M$,
in particular, for the geodesic flow corresponding to any
$G$-invariant Riemannian metric on
$M$. In general, a complete involutive subset of
$\{h\circ \mu^\mathrm{can}, h:\fg^*\to\bbR\}$
is not a complete involutive subset of the algebra
$C^\infty(T^*M)$. However, for the compact Lie group
$G$ the problem of constructing a complete involutive set of
real analytic functions on
$T^*(G/K)$ is reduced to the problem of finding a complete
involutive set of real analytic functions from the set
$A^G$ (see \cite[\S 2]{My01}, \cite[Lemma 3]{BJ01},
\cite{Pana03}). This is true also for the group
$\tilde G$ and the corresponding algebra
$A^{\tilde G}\subset C^\infty(T^*(\tilde G/\tilde K))$.

The algebra of functions
$A^G$ on $T^*(G/K)$, where, recall, $G/K$
is the adjoint orbit, contains some complete involutive subset
$\mathcal{F}$ of $A^G$
(\cite[Theorem 3.10]{MP04}, see also~\cite{BJ03},
\cite{BJ04}). The homogeneous space
$\tilde G/\tilde K$, as we remarked above, is a submanifold
of $G/K$ and therefore
$T(\tilde G/\tilde K)$ is a submanifold of
$T(G/K)$. Moreover,
$T(\tilde G/\tilde K)$ is a symplectic submanifold of
$T(G/K)$, where the symplectic structures on these spaces
are defined via isomorphisms
$T(\tilde G/\tilde K)\simeq T^*(\tilde G/\tilde K)$ and
$T(G/K)\simeq T^*(G/K)$ using a standard
$G$-invariant metric on $G/K$ and its restriction to
$\tilde G/\tilde K$ (see Proposition~\ref{pr.4}). The set
$\tilde{\mathcal{F}}=\{f|_{T(\tilde G/\tilde K)}, f\in\mathcal{F}\}$
of restrictions is an involutive subset of the algebra
$A^{\tilde G}$. This involutiveness of the functions from
$\tilde{\mathcal{F}}$ is a consequence of the fact that
$G/\tilde G$ is a symmetric space and follows easily from
results published in~\cite{MF78},~\cite{TF95}. The following
observation is crucial in our approach:\vskip4pt

\parbox{.9\textwidth}{
if the functions from the set $\mathcal{F}$ are independent at some
point of the symplectic submanifold
$T(\tilde G/\tilde K)\subset T(G/K)$, then the set
$\tilde{\mathcal{F}}$ is a complete
involutive subset of the algebra $A^{\tilde G}$.}
\vskip4pt

Therefore we investigate some open dense subset
$O$ of $T(G/K)$, where all functions from the set
$\mathcal{F}$ are independent (in the
paper~\cite{MP04} only the existence of such a set was
proved). We prove that
$O\cap T(\tilde G/\tilde K)\ne\emptyset$ if
$G$ is the unitary group $U(n)$ and $\sigma$
is its automorphism determined by the complex conjugation
(Theorem~\ref{th.16}).

The present paper is also motivated by the following
observation: in the above mentioned paper~\cite{DGJ09} the
first part of the proof of Theorem~4 needs an additional
argumentation (see more detailed comments in
Remark~\ref{re.17}). Nevertheless, the assertion of this
theorem is true as follows from our result.

Note that the existence of complete commutative
subalgebras of $A^G$ is an interesting and non-trivial problem.
In the review paper~\cite{BJ04a}
this problem was formulated as Conjecture 6.1
in terms of the so-called integrable pairs of Lie
groups $(G,K)$ (in the case of compact Lie groups).
The present paper proves this conjecture for a series of
important particular cases.

I thank the anonymous referees for a thorough reading of
the manuscript and a comprehensive list of suggestions that
helped me to improve the presentation.

\section{Some definitions, conventions, and notations}
\label{s.1}

All objects in this paper are real analytic,
$X$ stands for a connected manifold,
$\mathcal{E}(X)$ for the space of real analytic functions on
$X$.

We will say that some functions from the set
$\mathcal{E}(X)$ are independent if their differentials are
independent at each point of some open dense subset in
$X$. For any subset
${\mathcal F}\subset{\mathcal E}(X)$ denote by
$\ddim_x {\mathcal F}$ the maximal number of independent
functions from the set
${\mathcal F}$ at a point $x\in X$. Put
$\ddim{\mathcal F}\eqdef\max_{x\in X}\ddim_x{\mathcal F}$.

Let $\eta$ be a Poisson bi-vector on $X$ and let
$\mathcal{A}\subset\mathcal{E}(X)$
be a Poisson subalgebra of
$(\mathcal{E}(X),\eta)$, i.e. $\mathcal{A}$
is a real vector space closed under the Poisson bracket
$\{,\}:(f_1,f_2)\mapsto \eta(df_1,df_2)$ on $X$. Put
$(D\mathcal{A})_x\eqdef\{df_x: f\in \mathcal{A}\}\subset T^*_x X$
for any $x\in X$. Let $B_x$ denote the restriction of
$\eta_x$ to this subspace
$(D\mathcal{A})_x$. We say that a subset
$\mathcal{F}\subset\mathcal{A}$
is a {\it complete involutive subset of the algebra}
$(\mathcal{A},\eta)$ if at each point
$x$ of some open dense subset in $X$ the subspace
$V_x\eqdef\langle df_x: f\in
\mathcal{F}\rangle\subset(D\mathcal{A})_x$
is maximal isotropic with respect to the form
$B_x$, i.e. $B_x(V_x,V_x)=0$ and
$B_x(v,V_x)=0$ for $v\in(D\mathcal{A})_x$ implies
$v\in V_x$. In particular, any two functions
$f_1,f_2\in\mathcal{F}$ are in involution on
$X$, i.e. $\{f_1,f_2\}=0$.

\begin{definition}\label{de.1}
A pair $(\eta_1,\eta_2)$ of linearly independent bi-vector fields
(bi-vectors for short) on a manifold
$X$ is called Poisson if $\eta^t\eqdef t_1\eta_1+t_2\eta_2$
is a Poisson bi-vector for any
$t=(t_1,t_2)\in \bbR^2$, i.e. each bi-vector
$\eta^t$ determines on
$X$ a Poisson structure with the Poisson bracket
$\{,\}^{t}:(f_1,f_2)\mapsto \eta^t(df_1,df_2)$;
the whole family of Poisson bi-vectors
$\{\eta^t\}_{t\in\bbR^2}$ is called a {\it bi-Poisson
structure}.
\end{definition}

A bi-Poisson structure
$\{\eta^t\}$ (we will often skip the parameter space) can
be viewed as a two-dimensional vector space of Poisson
bi-vectors, the Poisson pair
$(\eta_1,\eta_2)$ as a basis in this space.

Suppose that a linear subspace
$\mathcal{A}\subset\mathcal{E}(X)$
is a Poisson subalgebra of
$(\mathcal{E}(X),\eta^t)$ for each
$t\in \bbR^2\setminus\{0\}$. Let
$B^t_x$ denote the restriction of
$\eta^t_x$ to the subspace $(D\mathcal{A})_x$,
$x\in X$.
\begin{definition}\label{de.2}
We say that the pair
$(\mathcal{A},\{\eta^t\})$ is {\it Kronecker} at a point
$x\in X$ for which
$\ddim_x \mathcal{A}=\ddim \mathcal{A}$ if the linear space
$\{B^t_x,t\in\bbC^2\}$ is two dimensional and
$\rank_\bbC B^t_x$ is constant with respect to
$(t_1,t_2)\in\bbC^2\setminus\{0\}$. We regard
$B^t_x,t\in\bbC^2$, as a complex bilinear form
$t_1 B^{(1,0)}_x+t_2 B^{(0,1)}_x$ on the complexification
$(D\mathcal{A})_x^\bbC\subset (T^*_x X)^\bbC$.
\end{definition}

The definitions above are motivated by the following assertion of
Bolsinov which is fundamental for our considerations.
\begin{proposition}\label{pr.3}~{\upshape\cite{Bol91}}
Let $B_1$ and $B_2$ be two linearly independent skew-symmetric
bilinear forms on a vector space
$V$. Suppose that the kernel of each form
$B^t=t_1B_1+t_2B_2$, $t\in\bbR^2$, is non-trivial, i.e.
$0<r\eqdef\min_{t\in\bbR^2}\dim\ker B^t$. Put
$T=\{t\in\bbR^2: \dim\ker B^t=r\}$. Then
\begin{mlist}
\item[{\upshape(1)}]
the subspace $L\eqdef\sum_{t\in T}\ker B^t$ is
isotropic with respect to any form $B^t$, $t\in\bbR^2$,
i.e. $B^t(L,L)=0$;
\item[{\upshape(2)}]
the space $L$ is maximal isotropic with respect to any
form $B^{t}$, $t\in T$ if and only if
$\dim_\bbC\ker B^t=r$ for all $t\in\bbC^2\setminus\{0\}$.
\end{mlist}
\end{proposition}

Consider a connected Riemannian manifold $(M,\mathbf{g})$ and its
connected Riemannian submanifold $(\tilde M,\tilde{\mathbf{g}})$, where
$\tilde{\mathbf{g}}=\mathbf{g}|_{\tilde M}$. The cotangent bundles $T^*M$
and $T^*\tilde M$ are symplectic manifolds with the canonical
symplectic structures $\Omega$ and $\tilde\Omega$
respectively. Using the metric $\mathbf{g}$ (resp. $\tilde{\mathbf{g}}$)
we can identify $T^*M$ with $TM$ (resp. $T^*\tilde M$ with $T\tilde
M$), denoting by $\varphi: TM\to T^*M$ (resp.
$\tilde\varphi: T\tilde M\to T^*\tilde M$ ) the corresponding
diffeomorphism. Let $p: TM\to M$ (resp. $\tilde p: T\tilde M\to \tilde M$)
be the natural projection and let
$\theta$ (resp. $\tilde\theta$) be the canonical 1-form on $T^*M$
(resp. on $T^*\tilde M$).
\begin{proposition}\label{pr.4}
The symplectic manifold $(T\tilde M,\tilde\varphi^*\tilde\Omega)$ is a
symplectic submanifold of $(TM,\varphi^*\Omega)$, i.e.
$\tilde\varphi^*\tilde\Omega=\varphi^*\Omega|_{T\tilde M}$. Moreover,
$\tilde\varphi^*\tilde\theta=\varphi^*\theta|_{T\tilde M}$.
\end{proposition}
\begin{proof}[Proof]
By the definition, $\theta_{x'}(Z')=x'(\pi_{*x'}Z')$,
where $\pi: T^*M\to M$ is the natural projection,
$x'\in T^*_qM$, $q=\pi(x')$, $Z'\in T_{x'}(T^*M)$.

If $x,y\in T_qM$, then $\varphi(x)\in T^*_qM$
and, by the definition, $\varphi(x)(y)=\mathbf{g}_q(x,y)$.
Putting $x'=\varphi(x)\in T^*_qM$ (to simplify the notation)
and taking into account that $\pi\circ\varphi=p$ we obtain
that for any $Z\in T_{x}TM$
\begin{equation*}
\begin{split}
(\varphi^*\theta)_{x}(Z)
&\eqdef
\theta_{\varphi(x)}(\varphi_{*x}Z)\eqdef
x'(\pi_{*x'}(\varphi_{*x}Z))\\
&=x'((\pi\circ\varphi)_{*x}Z)=
x'(p_{*x}Z)=
\mathbf{g}_q(x,p_{*x}Z).
\end{split}
\end{equation*}
Similarly, we obtain that
$$
(\tilde\varphi^*\tilde\theta)_{\tilde x}(\tilde Z)=
\tilde{\mathbf{g}}_{\tilde q}(\tilde x,\tilde p_{*\tilde x}\tilde Z)
\qquad\text{for any}\quad
\tilde q\in \tilde M, \
\tilde x\in T_{\tilde q}\tilde M,\
\tilde Z\in T_{\tilde x}T\tilde M.
$$
In other words,
$\tilde\varphi^*\tilde\theta=\varphi^*\theta|_{T\tilde M}$,
because $p|_{T\tilde M}=\tilde p$ and
$\mathbf{g}|_{\tilde M}=\tilde{\mathbf{g}}$.
Now to complete the proof it is sufficient to note that
$\Omega=d\theta$ and
$\tilde\Omega=d\tilde\theta$.
\qed
\end{proof}

\section{The integrability of geodesic flows}
\label{s.2}

In this section for any compact Lie algebra $\mathfrak{a}$
by $\fz(\mathfrak{a})$ we will denote its center and
by $\mathfrak{a}_s$ its maximal semisimple ideal, i.e.
$\mathfrak{a}=\fz(\mathfrak{a})\oplus\mathfrak{a}_s$; for any
real vector space or any Lie algebra $\mathfrak{a}$
by $\mathfrak{a}^\bbC$ we will denote its complexification.

\subsection{Commutator on $A^K_{\fm}$
induced by canonical Poisson structure on $T^*(G/K)$}
\label{ss.2.1}

Let $M=G/K$ be a homogeneous space of a compact connected Lie group
$G$  with the Lie algebra
$\fg$. There exists a faithful representation $\chi$ of
$\fg$ such that its associated bilinear form
$\Phi_\chi$ is negative-definite on $\fg$ (if $\fg$
is semi-simple we can take the Killing form associated with
the adjoint representation of $\fg$).
Let $\fm=\fk^\bot$ be the orthogonal complement to $\fk$ with
respect to $\Phi_\chi$. Then
\begin{equation}\label{eq.1}
\fg=\fm\oplus\fk,\qquad [\fk,\fm]\subset\fm.
\end{equation}

The form $\langle,\rangle=-\Phi_\chi$ determines a $G$-invariant
metric on $G/K$. This metric identifies the cotangent
bundle $T^*(G/K)$ and the tangent bundle $T(G/K)$.
Thus we can also talk on the canonical 2-form
$\Omega$ on the manifold $T(G/K)$. The symplectic form $\Omega$ is
$G$-invariant with respect to the natural action of $G$
on $T(G/K)$ (extension of the action of $G$ on $G/K$).

We can identify the tangent space
$T_o(G/K)$ at the point $o=p(e)$ with the space $\fm$ by
means of the canonical projection $p:G\to G/K$. Let
$A^G$ (resp. $A^K_{\fm}$) be the set of all
$G$-invariant (resp. $\Ad(K)$-invariant)
functions on $T(G/K)$ (resp. on $\fm$). There is a
one-to-one correspondence between $G$-orbits in $T(G/K)$
and $\Ad(K)$-orbits in $\fm$. Thus we can naturally identify the
spaces of functions $A^G$ and $A^K_{\fm}$.
For any smooth function
$f$ on $\fm$ write $\grad_\fm f$
for the vector field on $\fm$ such that
$$
df_x(y)=\langle \grad_\fm f(x), y \rangle
 \qquad\text{for all}\qquad
y\in\fm.
$$

The Poisson bracket of two functions
$f_1,f_2$ from the set $A^K_\fm=A^G$
with respect to the canonical Poisson structure
$\eta^\mathrm{can}$ (determined by the canonical 2-form
$\Omega$) has the form~\cite[Lemma 3.1]{MP04} :
\begin{equation}\label{eq.2}
\{f_1,f_2\}^\mathrm{can}(x)= -\langle x,\
[\grad_\fm f_1(x),\grad_\fm f_2(x)]\rangle,
\quad
x\in\fm.
\end{equation}

Now, let us consider an important for our considerations subspace of
$\fm$. For any $x\in\fm$ define the subspace
$\fm(x)\subset\fm$ putting
\begin{equation}\label{eq.3}
\fm(x)\eqdef\{y\in\fm: [x,y]\in\fm\}=
\{y\in\fm:\langle y,\ad x(\fk)\rangle=0\},
\end{equation}
in particular,
$$
\ad x\bigl(\fm(x)\bigr)\subset\fm
\quad\text{and}\quad \fm(x)\oplus\ad x(\fk)=\fm.
$$

For any element $x\in\fg$ denote by $\fg^x$
its centralizer in $\fg$, i.e. the set of all $z\in\fg$
satisfying $[x,z]=0$. Put $\fk^x=\fg^x\cap \fk$.
Consider in $\fm$ a nonempty Zariski open subset:
\begin{equation}\label{eq.4}
R(\fm)
=\{x\in\fm :\dim\fg^x=q(\fm), \dim\fk^x=p(\fm)\},
\end{equation}
where $q(\fm)$ (resp.
$p(\fm)$) is the minimum of dimensions of the spaces
$\fg^y$ (resp.$\fk^y$) over all $y\in\fm$.

Let $(\cdot)_\fm$ be the projection of $\fg$ into
$\fm$ along $\fk$. For each
$x\in R(\fm)$ the spaces $\fm(x)$ and
$(\fg^x)_\fm\subset\fm(x)$
have the same dimensions:
$$
\dim\fm(x)= \dim\fm-\dim\ad x(\fk)
\quad\text{and}\quad
\dim (\fg^x)_\fm=\dim\fg^x-\dim\fk^x.
$$
Moreover, for each
$x\in R(\fm)$ the maximal semi-simple ideal
$(\fg^x)_s=[\fg^x,\fg^x]$ of
$\fg^x$ is contained in the algebra $\fk^x$, i.e.
\begin{equation}\label{eq.5}
[\fg^x,\fg^x]=[\fk^x,\fk^x],\quad
\text{if}\ x\in R(\fm),
\end{equation}
(see~\cite[Prop.10]{My01} or~\cite{Mi82}).
Therefore, $\dim(\fg^x/\fk^x)=\rank\fg-\rank\fk^x$, i.e.
\begin{equation}\label{eq.6}
\dim\fg^x=\rank\fg+(\dim\fk^x-\rank\fk^x),
\qquad \text{if}\ x\in R(\fm).
\end{equation}

It is clear that $\grad_\fm f(x)\in\fm(x)$
for any $f\in A^K_\fm$.
Moreover, since the Lie group
$K$ is compact, for each
$x$ from some nonempty Zariski open subset of
$\fm$ the space $\fm(x)$ is generated by the vectors
$\grad_\fm f(x)$, $f\in A^K_\fm$.
For any $x\in R(\fm)$ the kernel of the 2-form
$B^0_x:\fm(x)\times\fm(x)\to\bbR$,
$(y_1,y_2)\mapsto -\langle x,[y_1,y_2]\rangle$,
associated with the Poisson structure $\eta^\mathrm{can}$,
is the space $(\fg^x)_\fm$:
\begin{equation*}
\begin{split}
\{y\in\fm(x)&: \langle x,[y,\fm(x)]\rangle=0\}
=\{y\in\fm(x): \langle [x,y],\fm(x)\rangle=0\} \\
&=\{y\in\fm(x): [x,y]\in\ad x(\fk)\}=(\fg^x)_\fm\cap \fm(x)=(\fg^x)_\fm.
\end{split}
\end{equation*}
Hence the number $\frac12(r(\fm)+\dim\fm(x))$, where
\begin{equation}\label{eq.7}
r(\fm)=\dim (\fg^x)_\fm=\dim\fg^x-\dim\fk^x=
q(\fm)-p(\fm),\quad
x\in R(\fm),
\end{equation}
is the maximal number of functions in involution from the set
$A^K_\fm$ functionally independent at the point $x$.

For an arbitrary element $x\in R(\fm)$, we have~\cite[Prop.9]{My01}
\begin{equation}\label{eq.8}
[\fm(x),\fk^x]=0.
\end{equation}

The compact Lie algebra
$\fk$ is the direct sum $\fk=\fz(\fk)\oplus\fk_s$ of the center and
of the maximal semisimple ideal. The center $\fz(\fk)$ of
$\fk$ will be denoted simply by $\fz$ for short. Then we
have the following orthogonal splittings with respect to the
invariant form $\langle\cdot,\cdot \rangle$:
$\fk=\fz\oplus\fk_s$, $\fg=(\fm\oplus\fz)\oplus\fk_s$.

Consider the set $R(\fm\oplus\fz)$ determined
by~(\ref{eq.4}) for the pair $(\fg,\fk_s)$.
Then for any $x+z\in\fm\oplus\fz$ such that
$x\in R(\fm)$, $z\in\fz$ and $x+z\in R(\fm\oplus\fz)$,
we have $\fk_s^{x+z}=\fk_s^x$ because $[\fz,\fk]=0$. However,
$\fk_s$ contains the maximal semi-simple ideal
of the compact Lie algebra $\fk^x$, and therefore
$(\dim\fk^x-\rank\fk^x)=(\dim\fk_s^x-\rank\fk_s^x)$.
Thus by~\eqref{eq.6}
\begin{equation}\label{eq.9}
\begin{split}
q(\fm)=\dim\fg^x
&=\rank\fg+(\dim\fk^x-\rank\fk^x)
=\rank\fg+(\dim\fk_s^x-\rank\fk_s^x)\\
&=\rank\fg+(\dim\fk_s^{x+z}-\rank\fk_s^{x+z})
=\dim\fg^{x+z}=q(\fm\oplus\fz).
\end{split}
\end{equation}

The following proposition will be often used in
subsections~\ref{ss.2.2} and~\ref{ss.2.3}.
\begin{proposition}\label{pr.5}
Suppose that $\fk= \fg^a$ for some $a\in\fg$ and
$x_0\in R(\fm)$. Let $\fg_0$ and $\fk_0$
be the centralizers of the algebra $\fk^{x_0}$ in $\fg$ and $\fk$
respectively. Let $\fm_0=\fg_0\cap\fm$. Then
\begin{mlist}
\item[{\upshape(1)}]
for any $x\in\fm_0\cap R(\fm)$ {\upshape(}this set contains
$x_0${\upshape)} we have
$\fm_0(x)=\fm(x);$
\item[{\upshape(2)}]
for any $x\in\fm_0\cap R(\fm)$ we have $\fk^x=\fk^{x_0}$
and $x\in R(\fm_0);$
\item[{\upshape(3)}]
$a\in\fg_0$, $\fg^{x_0}=\fg_0^{x_0}\oplus(\fk^{x_0})_s$
and  the direct sum $\fg_0\oplus(\fk^{x_0})_s$ is
a subalgebra of maximal rank of $\fg;$
\item[{\upshape(4)}]
each element $x$ of the set $R(\fm_0)$
is a regular element of the Lie
algebra $\fg_0$ and its centralizer $\fk_0^x$ in $\fk_0$ is
the center $\fz(\fg_0)$ of $\fg_0;$
\item[{\upshape(5)}]
$r(\fm)=r(\fm_0)=\rank\fg_0-\dim\fz(\fg_0).$
\end{mlist}
\end{proposition}
\begin{proof}[Proof]
Items $(1)$ and $(2)$ follows immediately from
Proposition~2.3~\cite{MP04}. It is clear that
$a\in\fg_0$ because $a\in\fz(\fk)$. By~(\ref{eq.5}),
$(\fg^{x_0})_s=(\fk^{x_0})_s\subset\fk$, i.e.
$\fg^{x_0}=\fz(\fg^{x_0})\oplus(\fk^{x_0})_s$ and
$\fk^{x_0}=\fz(\fk^{x_0})\oplus(\fk^{x_0})_s$.
Since by definition $\fk^{x_0}\subset\fg^{x_0}$, we obtain that
$[\fz(\fg^{x_0}),\fk^{x_0}]=0$, i.e.
$\fz(\fg^{x_0})\subset\fg_0$, and consequently,
$\fg^{x_0}_0\eqdef\fg^{x_0}\cap\fg_0=\fz(\fg^{x_0})$.
In other words, $x_0$ is a regular element of the Lie algebra
$\fg_0$ and $\fg^{x_0}=\fg_0^{x_0}\oplus(\fk^{x_0})_s$.
Since the Lie algebra
$(\fk^{x_0})_s$ is semisimple,
$\fg_0\cap(\fk^{x_0})_s=0$. The subalgebra
$\fg_0\oplus(\fk^{x_0})_s\subset\fg$
is a subalgebra of maximal rank because it contains the centralizer
$\fg^{x_0}$.

Each element $x$ of the set $R(\fm_0)$ is also regular in
$\fg_0$ because $x_0\in R(\fm_0)$ by
$(2)$ and by definition of the set
$R(\fm_0)$ (see~(\ref{eq.4})),
$\dim\fg^{x}_0=\dim\fg^{x_0}_0$ ($=\rank\fg_0$).

By the definition, $\fz(\fk^{x_0})=\fk^{x_0}\cap\fg_0$ and
$\fz(\fk^{x_0})\subset\fz(\fg_0)$. Since
$\fk^{x_0}_0\eqdef\fk^{x_0}\cap\fg_0$, we obtain that
$\fk^{x_0}_0=\fz(\fk^{x_0})$ and
$\fk^{x_0}_0\subset\fz(\fg_0)$. But
$a,x_0\in\fg_0$. Hence these elements commute with the
center $\fz(\fg_0)$ of $\fg_0$ and, consequently,
$\fz(\fg_0)\subset\fg^a=\fk$ and
$\fz(\fg_0)\subset\fg^{x_0}$. In other words,
$\fz(\fg_0)\subset\fk\cap\fg^{x_0}\cap\fg_0\eqdef\fk^{x_0}_0$,
i.e. $\fk^{x_0}_0=\fz(\fg_0)$. Since
$\dim\fk_0^{x_0}=\dim\fk_0^{x}$ for $x\in R(\fm_0)$ and
$\fz(\fg_0)\subset(\fg^{x}_0\cap\fk)\eqdef\fk^{x}_0$, item
$(4)$ is proved. To complete the proof of the proposition it
is sufficient to remark that
$x_0\in R(\fm)\cap R(\fm_0)$,
$$
r(\fm)
=\dim\fg^{x_0}-\dim\fk^{x_0}
=\dim\fz(\fg^{x_0})-\dim\fz(\fk^{x_0})
=\dim\fg^{x_0}_0-\dim\fk^{x_0}_0
=r(\fm_0)
$$
and $x_0$ is a regular element of $\fg_0$, i.e.
$\dim\fg^{x_0}_0=\rank\fg_0$. \qed
\end{proof}

\subsection{ The bi-Poisson structure
$\{\eta^t(\fg,a)\}$: exact
formulas and involutive sets of functions}\label{ss.2.2}

Consider the adjoint action
$\Ad$ of $G$ on the Lie algebra
$\fg$. Suppose now in addition that the Lie subgroup
$K$ is an isotropy group of some element
$a\in\fg$, i.e. $K=\{g\in G: \Ad g(a)=a\}$ and
$\fk=\fg^a$. Moreover, by invariance of the form
$\langle\cdot ,\cdot \rangle$, $\fk=\fz(\fk)\oplus\fk_s$,
$a\in\fz(\fk)$, and $\ad a(\fm)\subset\fm$.

Using the invariant form
$\langle\cdot,\cdot\rangle$ on the Lie algebra
$\fg$, we identify the dual space $\fg^*$ and
$\fg$. So the orbit ${\mathcal O}=G/K$
is a symplectic manifold with the Kirillov-Kostant-Souriau
symplectic structure
$\omega_{\mathcal O}$. By the definition, the form
$\omega_{\mathcal O}$ is $G$-invariant and at the point
$a\in{\mathcal O}$ we have
$$
\omega_{\mathcal O}(a)([a,\xi_1],[a,\xi_2])=-\langle
a,[\xi_1,\xi_2]\rangle, \quad \forall\xi_1,\xi_2\in\fg,
$$
where we consider the vectors
$[a,\xi_1],[a,\xi_2]\in T_a\fg=\fg$
as tangent vectors to the orbit
${\mathcal O}\subset\fg$ at the point
$a\in{\mathcal O}$. Let $\tau: T{\mathcal O}\to{\mathcal O}$
be the natural projection. Using the closed 2-form
$\tau^*\omega_{\mathcal O}$ on $T{\mathcal O}$ (the lift of
$\omega_{\mathcal O}$) we construct a bi-Poisson structure
on $T{\mathcal O}$.

Consider on $T{\mathcal O}$ two symplectic forms:
$\omega_1=\Omega$ and
$\omega_2=\Omega+\tau^*\omega_{\mathcal O}$. Write
$\eta_1=\omega_1^{-1}$, $\eta_2=\omega_2^{-1}$
for the standard Poisson bi-vectors associated with these forms, i.e.
$\eta_i(df_1,df_2)=-\omega_i(\xi_{f_1},\xi_{f_2})$, where
$\xi_{f_j}$, are the Hamiltonian vector fields of the
functions $f_j$ ($df_j=-\omega_i(\xi_{f_j},\cdot)$),
$i,j=1,2$. Then the family
$\{\eta^t(\fg,a)=\eta^t=t_1\eta_1+t_2\eta_2\}$,
$t_1,t_2\in\bbR$, is a bi-Poisson structure on
$T\mathcal{O}$~\cite[Prop.1.6]{MP04}. Putting
$t_2=\lambda$, $t_1=1-\lambda$,
$\lambda\in\bbR$ or $t_1=-1$, $t_2=1$
we exclude considering of proportional bi-vectors. The
corresponding bi-vectors are denoted by
$\eta^\lambda$, $\lambda\in\bbR$, and
$\eta^\mathrm{si}$ (the singular bi-vector). The Poisson
bracket of two functions
$f_1,f_2$ from the set
$A^K_\fm=A^G\subset {\mathcal E}(T{\mathcal O})$
with respect to the Poisson structure
$\eta^\lambda$, $\lambda\in\bbR$ or
$\eta^\mathrm{si}$ has the form~\cite[Lemma 3.1]{MP04}:
\begin{equation*}
\begin{split}
\{f_1,f_2\}^\lambda(x)
&= -\langle x+\lambda a,\ [\grad_\fm f_1(x),\grad_\fm f_2(x)]\rangle,\\
\{f_1,f_2\}^\mathrm{si}(x)
&= -\langle a,\ [\grad_\fm f_1(x),\grad_\fm f_2(x)]\rangle.
\end{split}
\end{equation*}
Remark that the structure $\eta^0$ ($\lambda=0$)
is the canonical Poisson structure (see~\eqref{eq.2}).
Since the bi-Poisson structure
$\{\eta^t\}$ is $G$-invariant
it is sufficient to investigate points in $R(\fm)\subset\fm=T_o(G/K)$, where the
pair $(A^K_\fm,\eta^t)$ is Kronecker.

The bi-Poisson structure
$\{\eta^t\}$ determines at a point
$x\in\fm=T_o(G/K)$ the bilinear forms
$B^t_x: (DA^G)_x\times (DA^G)_x\to\bbR$, where, recall,
$B^t_x$ is the restriction $\eta^t|_{(DA^G)_x}$. Let
$x$ be an element of
$R(\fm)\subset T_o(G/K)$. Since we identified the spaces
$A^G$ and $A^K_\fm$, $B^t_x$ determines the following
complex-valued bilinear forms (which we denote also by
$B^t_x$, $B^\lambda_x$ and $B^\mathrm{si}_x$ for short) on
$\fm(x)\times\fm(x)$~\cite[(3.11)]{MP04}:
\begin{equation}\label{eq.10}
\begin{split}
B^t_x
&:(y_1,y_2)\mapsto -\langle (t_1+t_2)x+t_2a,\
[y_1,y_2]\rangle,\ t_1,t_2\in\bbC,\\
B^\lambda_x
&:(y_1,y_2)\mapsto -\langle x+\lambda a,\ [y_1,
y_2]\rangle,\quad \lambda\in\bbC, \\
B^\mathrm{si}_x
&:(y_1,y_2)\mapsto -\langle a,\ [y_1,y_2]\rangle.
\end{split}
\end{equation}
Let $\fm^\bbC(x)$ be the complexification of the space
$\fm(x)$, $x\in R(\fm)$. It is easy to see that the kernel
of the form $B^\lambda_x$ in $\fm^\bbC(x)$ is the subspace of
$\fm^\bbC(x)$ given by
\begin{equation*}
\begin{split}
\ker B^\lambda_x
&=\{y\in\fm^\bbC(x): [x+\lambda a,\ y]\in\ad x(\fk^\bbC)\} \\
&=\{y\in\fm^\bbC(x): [x+\lambda a,\ y]\in\ad (x+\lambda a)(\fk^\bbC)\}
\end{split}
\end{equation*}
because $\ad x(\fk^\bbC)=(\fm^\bbC(x))^\bot$ in $\fm^\bbC$
and $[a,\fk^\bbC]=0$. Thus
$$
\ker B^\lambda_x=((\fg^\bbC)^{x+\lambda a})_{\fm^\bbC}\cap\fm^\bbC(x),
$$
where $(\cdot)_{\fm^\bbC}$ denotes the projection onto
$\fm^\bbC$ along $\fk^\bbC$. However,
$((\fg^\bbC)^{x+\lambda a})_{\fm^\bbC}\subset\fm^\bbC(x)$
because $[a,\fm^\bbC]\subset\fm^\bbC$, $[a,\fk^\bbC]=0$,
$\ad x(\fk^\bbC)\subset\fm^\bbC$ and by~\eqref{eq.3}
$y\in\fm^\bbC$ is an element of $\fm^\bbC(x)$ if and only if
$[x,y]\in\fm^\bbC$. Thus,
$$
\ker B^\lambda_x =((\fg^\bbC)^{x+\lambda a})_{\fm^\bbC},
\ \lambda\in\bbC.
$$
 In particular, for
$\lambda=0$ (for the canonical Poisson structure on
$T(G/K)$),
$$
\ker B^0_x=((\fg^\bbC)^x)_{\fm^\bbC}=((\fg^x)_\fm)^\bbC.
$$
Since $x\in R(\fm)$, the (real) dimension of the space
$(\fg^x)_\fm$ is equal to the constant
$r(\fm)=q(\fm)-p(\fm)$.
Therefore a maximal isotropic
subspace of the space $\fm(x)$ with respect to the form
$B^0_x$ is of dimension $\frac12\bigl(r(\fm)+\dim\fm(x)\bigr)$.
It is clear that
$$
\ker B^\mathrm{si}_x
=\{y\in\fm^\bbC(x): [a,\ y]\in\ad x(\fk^\bbC)\}
=\fm^\bbC(x)\cap \bigl(\ad_a^{-1}\ad x(\fk^\bbC)\bigr),
$$
where $\ad_a^{-1}\eqdef (\ad a|_\fm)^{-1}$.
As an immediate consequence of Proposition~\ref{pr.3} we
obtain
\begin{proposition}\label{pr.6}
The pair
$(A^G,\eta^t(\fg,a))$ is Kronecker at a point
$x\in R(\fm)$ if and only if

$(1)$ $\dim_\bbC ((\fg^\bbC)^{x+\lambda a})_{\fm^\bbC}=r(\fm)$
for each $\lambda\in\bbC${\rm ;}

$(2)$
$\dim_\bbC\ker B^\mathrm{si}_x=r(\fm)$.
\end{proposition}

Let $O^\mathrm{Kr}(\fm)$ be a set of all points
$x\in R(\fm)$ for which conditions $(1)$ and
$(2)$ of Proposition~\ref{pr.6} hold. By
Proposition~\ref{pr.3} this set is a Zariski open (possibly
empty) subset of $\fm$.

Denote by $I(\fg)$ the space of all $\Ad(G)$-invariant
polynomials on $\fg$. If $h\in I(\fg)$ then it is clear
that the function $h^{\lambda}:\fg\to\bbR$,
$h^{\lambda}(y)=h(y+\lambda a)$, $\lambda\in\bbR$, is
$\Ad(K)$-invariant on $\fg$.
Therefore the set
$\mathcal{F}(\fg,\fm)=\{h^\lambda|_\fm, h\in I(\fg),\ \lambda\in\bbR\}$
is a subset of $A^{K}_{\fm}=A^{G}$
(of $G$-invariant functions on $T(G/K)$).
The following assertion was proved in~\cite{MP04}
(see Proposition~3.6, Lemma~3.3 and Theorem~3.9).
\begin{theorem}\label{th.7}{\upshape\cite{MP04}}
The Zariski open subset
$O^\mathrm{Kr}(\fm)\subset R(\fm)$ of
$\fm$ is nonempty. The set of functions
$\mathcal{F}(\fg,\fm)$ is a complete involutive subset of the Poisson
algebra $(A^G,\eta^0)$. Moreover, for each point
$x\in O^\mathrm{Kr}(\fm)\subset R(\fm)$ the space
$\{\grad_\fm f(x), f\in\mathcal{F}(\fg,\fm)\}\subset\fm(x)$
is a maximal isotropic subspace of
$\fm(x)$.
\end{theorem}
By Theorem~\ref{th.7} the set
\begin{equation}\label{eq.11}
\begin{split}
Q_a(\fm)
&=\{x\in R(\fm): \dim_\bbC\ker B^\mathrm{si}_x=r(\fm)\} \\
&=\{x\in R(\fm): \dim\left( \fm(x)\cap
\bigl(\ad_a^{-1}\ad x(\fk)\bigr)\right)=r(\fm)\}
\end{split}
\end{equation}
contains the dense open subset
$O^\mathrm{Kr}(\fm)\subset\fm$, in particular,
\begin{equation}\label{eq.12}
r(\fm)=\min_{x\in R(\fm)}\dim \Bigl(\fm(x)\cap
\bigl(\ad_a^{-1}\ad x(\fk)\bigr)\Bigl).
\end{equation}
Since $\fm(x)=(\ad x(\fk))^\bot$ in
$\fm$ and $\dim\fm(x)$ is constant for $x\in R(\fm)$,
the set $Q_a(\fm)\subset R(\fm)$ is a nonempty Zariski open subset of
$\fm$. Put
\begin{equation}\label{eq.13}
M_a(\fm)=\{x\in \fm:
\dim_\bbC (\fg^\bbC)^{x+\lambda a}=q(\fm)
\text{ for each } \lambda\in\bbC\}.
\end{equation}
Then
\begin{equation}\label{eq.14}
O^\mathrm{Kr}(\fm)=Q_a(\fm)\cap M_a(\fm)
\end{equation}
because $Q_a(\fm)\subset R(\fm)$ and
$(\fk^\bbC)^{x+\lambda a}=(\fk^\bbC)^{x}$ for each
$\lambda\in\bbC$.
\begin{lemma}\label{le.8}
The set $M_a(\fm)$ is a Zariski open subset of $\fm$.
\end{lemma}
\begin{proof}[Proof]
The assertion of the lemma in more general form was used
in the paper~\cite[\S 3.4]{BJ04} (but without any proof).
For completeness we will prove this lemma here.

Let $V=\fm^\bbC\oplus\bbC a\subset\fg^\bbC$.
Since by~(\ref{eq.9}),
$\min\limits_{x\in\fm^\bbC,
\lambda\in\bbC}\dim_\bbC(\fg^\bbC)^{x+\lambda a}=q(\fm)$,
the set
$$S=\{y\in V: \dim_\bbC(\fg^\bbC)^y>q(\fm)\}$$
is a Zariski closed subset of $V$.
If $y\in S$, then any non-zero (complex) scalar multiple of
$y$ is also an element of $S$. So that one can consider a closed subvariety
$\mathbb{P} S$ in the projectivization $\mathbb{P} V$ of $V$.

Suppose that $M_a(\fm)$ is nonempty and
$x_0\in M_a(\fm)$. Then by definition
$(x_0+\bbC a)\cap S=\emptyset$.
In other words, the intersection of the projective line
$\mathbb{P}\langle x_0,a\rangle$ with the set
$\mathbb{P} S$ is empty. Clearly, any nearby line
$\mathbb{P}\langle x,a\rangle$,
$x\in\fm^\bbC$ will also not intersect
$\mathbb{P} S$. The set $M_a(\fm^\bbC)$ of all such
$x\in\fm^\bbC$ is Zariski open in
$\fm^\bbC$. Taking into account that
$M_a(\fm)=M_a(\fm^\bbC)\cap\fm$ we complete the proof. \qed
\end{proof}

\subsection{Adjoint orbits and involutive automorphisms}
\label{ss.2.3}

Let $\sigma$ be an involutive automorphism of
$\fg$ and let
$\fg=\tilde\fg\oplus\fg'$ be the decomposition of
$\fg$ into the eigenspaces of
$\sigma$ for the eigenvalues $+1$ and
$-1$ respectively:
$$
[\tilde\fg,\tilde\fg]\subset\tilde\fg,
\qquad
[\fg',\fg']\subset\tilde\fg,
\qquad
[\tilde\fg,\fg']\subset\fg'.
$$
Denote by $\tilde G$ the closed connected subgroup of
$G$ with the Lie algebra
$\tilde\fg$. Fix some element $a\in\fg'$
($\sigma(a)=-a$) and consider the orbit
$\tilde{\mathcal O}=\Ad(\tilde G)(a)=\tilde G/\tilde K$ in
$\fg$ through this element $a$. It is clear that
$\tilde{\mathcal O}$ is a submanifold
($\tilde G$-suborbit) of the $G$-orbit
${\mathcal O}=G/K$ of $a$ and
$\tilde K=\tilde G\cap K$. The form
$\langle\cdot,\cdot\rangle$ determines a
$\tilde G$-invariant metric on
$\tilde{\mathcal O}$. This metric identifies the cotangent
bundle $T^*\tilde{\mathcal O}$ and the tangent bundle
$T\tilde{\mathcal O}$. Thus we can also talk on the
canonical 2-form $\tilde\Omega$ on the manifold
$T\tilde{\mathcal O}$.

Since $\sigma(a)=-a$, the algebra $\fk=\fg^a$ is
$\sigma$-invariant. Suppose that the form
$\langle\cdot,\cdot\rangle=-\Phi_\chi$ is also
$\sigma$-invariant (if $\fg$ is semi-simple its Killing
form is invariant with respect to an arbitrary automorphism of
$\fg$). Then $\sigma(\fm)=\fm$ and we have in addition
to~\eqref{eq.1} the following orthogonal decompositions of algebras
$\fg,\tilde\fg,\fk$ with respect to the form
$\langle\cdot ,\cdot \rangle$:
\begin{equation}\label{eq.15}
\fg=\tilde\fk\oplus\fk'\oplus\tilde\fm\oplus\fm',
\quad
\tilde\fg=\tilde\fk\oplus\tilde\fm,
\quad
[\tilde\fk,\tilde\fm]\subset\tilde\fm,
\quad
[\tilde\fk,\fm']\subset\fm',
\quad
\fk=\tilde\fk\oplus\fk',
\end{equation}
where $\tilde\fk,\tilde\fm$ are subspaces of
$\tilde\fg$, $\fk',\fm'$ are subspaces of
$\fg'$. In particular, $(\fk,\tilde\fk)$
is a symmetric pair of compact Lie algebras, i.e.
$\tilde\fk$ is the fixed point set of the involutive
automorphism $\sigma|_\fk$.

Since $\ker\ad a=\fk$ and $\fm=\fk^\bot$ in $\fg$, then
$\ad a(\fm)=\fm$ and the operator
$\ad a|_\fm:\fm\to\fm$ is invertible. Moreover, for
$\fm'\subset\fm$ and $\tilde\fm\subset\fm$ we have
$$
\ad a(\fm')\subset[\fg',\fm']
\cap\fm\subset\tilde\fg\cap\fm\subset\tilde\fm,
\quad\text{and}\quad
\ad a(\tilde\fm)\subset[\fg',\tilde\fm]
\cap\fm\subset\fg'\cap\fm\subset\fm',
$$
and therefore $\dim\tilde\fm=\dim\fm'$. Since
$\ker(\ad a|_\fm)=0$, we have
\begin{equation}\label{eq.16}
\ad a(\fm')=\tilde\fm,
\quad\text{and}\quad
\ad a(\tilde\fm)=\fm'.
\end{equation}

Similarly as in the case of the pair $(\fg,\fk)$
for any $x\in\tilde\fm$ define the subspace
$\tilde\fm(x)\subset\fm$ putting
$$
\tilde\fm(x)\eqdef\{y\in\tilde\fm: [x,y]\in\tilde\fm\}=
\{y\in\tilde\fm:\langle y,\ad x(\tilde\fk)\rangle=0\}.
$$
Then, in particular, $\ad x\bigl(\tilde\fm(x)\bigr)\subset\tilde\fm$
and $\tilde\fm(x)\oplus\ad x(\tilde\fk)=\tilde\fm$.

Let $A^{\tilde G}$ (resp.
$A^{\tilde K}_{\tilde\fm}$) be the set of all
$\tilde G$-invariant (resp.
$\Ad(\tilde K)$-invariant) functions on
$T(\tilde G/\tilde K)$ (resp. on
$\tilde\fm$). The Poisson bracket of two functions
$\tilde f_1,\tilde f_2\in A^{\tilde K}_{\tilde\fm}$
with respect to the canonical Poisson structure
$\tilde\eta^\mathrm{can}$ (determined by the canonical
2-form $\tilde\Omega$ on
$T\tilde{\mathcal O}$) has the form (see~(\ref{eq.2})) :
\begin{equation}\label{eq.17}
\{\tilde f_1,\tilde f_2\}^\mathrm{can}(x)
= -\langle x,\ [\grad_{\tilde\fm} \tilde f_1(x),
\grad_{\tilde\fm} \tilde f_2(x)]\rangle,
\quad
x\in\tilde\fm.
\end{equation}
Note that by $\Ad(\tilde K)$-invariance of the function
$\tilde f_k$ the vector $\grad_{\tilde\fm} \tilde f_k(x)$,
$k=1,2$, belongs to the subspace
$\tilde\fm(x)\subset\tilde\fm$.

Since $\sigma(a)=-a$, the center
$\fz=\fz(\fk)$ of the reductive Lie algebra
$\fk=\fg^a$ is $\sigma$-invariant, i.e.
$\fz=\tilde\fz\oplus\fz'$, where
$\tilde\fz=\fz\cap\tilde\fg$,
$\fz'=\fz\cap\fg'$. It is clear that
$a\in\fz'$. Then for each element
$b\in\fz'$ we can consider the endomorphism
$\varphi_{a,b}:\fg\to\fg$ on $\fg$ putting
$\varphi_{a,b}(x)=\ad^{-1}_a([b,x])$ for
$x\in\fm$ and $\varphi_{a,b}(z)=z$ for
$z\in\fk$, where, recall,
$\ad_a^{-1}\eqdef(\ad a|_\fm)^{-1}$. Remark that
$\varphi_{a,b}(\fm)\subset\fm$ because
$[\fk,\fm]\subset\fm$. Moreover,
$\varphi_{a,b}(\tilde\fm)\subset\tilde\fm$ and
$\varphi_{a,b}(\fm')\subset\fm'$ because
$a,b\in\fg'$ (see also~(\ref{eq.16})).

It is clear that the endomorphism
$\varphi_{a,b}$ is symmetric and the group
$\Ad(K)$ commutes elementwise with
$\varphi_{a,b}$ on $\fm$. Therefore, the operator
$\varphi_{a,b}|_{\tilde\fm}$ is also symmetric and the group
$\Ad(\tilde K)$ commutes elementwise with
$\varphi_{a,b}$ on $\tilde\fm$
($K$ is connected). Therefore the function
$\tilde H_{a,b}(x)=\frac12\langle x,\varphi_{a,b}(x)\rangle$,
$x\in\tilde\fm$, is $\Ad(\tilde K)$-invariant. Then
$\tilde H_{a,b}$ (as a function on
$T(\tilde G/\tilde K)$ from the set
$A^{\tilde G}=A^{\tilde K}_{\tilde\fm}$)
is a Hamiltonian function of the geodesic flow of some
pseudo-Riemannian metric on
$\tilde G/\tilde K$ if
$\varphi_{a,b}|_{\tilde\fm}$ is non-degenerate.

Consider the space
$I(\fg)$ of all $\Ad(G)$-invariant polynomials on
$\fg$. As we remarked in the previous subsection, for each
$h\in I(\fg)$ the function $h^{\lambda}:\fg\to\bbR$,
$h^{\lambda}(y)=h(y+\lambda a)$, $\lambda\in\bbR$, is an
$\Ad(K)$-invariant function on $\fg$. Therefore the set
$\mathcal{F}(\fg,\fm)=\{h^\lambda|_\fm, h\in I(\fg),\,\lambda\in\bbR\}$
is a subset of $A^{K}_{\fm}=A^{G}$ (of $G$-invariant function on
$T(G/K)$) and the set
$\mathcal{F}(\fg,\tilde\fm)=\{h^\lambda|_{\tilde\fm}, h\in
I(\fg),\, \lambda\in\bbR\}$
is a subset of $A^{\tilde K}_{\tilde\fm}=A^{\tilde G}$ (of
$\tilde G$-invariant function on
$T(\tilde G/\tilde K)$). Put $H^\lambda=h^\lambda|_\fm$ and
$\tilde H^\lambda=h^\lambda|_{\tilde\fm}$.
The following lemma follows easily from the results
of~\cite{MF78} (see also~\cite[Ch.6,16.Lemma]{TF95}
or~\cite[sec.3]{DGJ09}).
\begin{lemma}\label{le.9}~{\upshape\cite{DGJ09}}
For any functions $h_1,h_2$, $h\in I(\fg)$ and arbitrary
parameters $\lambda_1,\lambda_2,\lambda\in\bbR$ we have
$\{\tilde H_1^{\lambda_1},\tilde H_2^{\lambda_2}\}^\mathrm{can}=0$ and
$\{\tilde H^{\lambda},\tilde H_{a,b}\}^\mathrm{can}=0$.
\end{lemma}
\begin{proof}[Proof]
Mainly to fix notations we will prove this lemma here. Since
$\sigma$ is an automorphism of $\fg$ and
$\Ad(G)$ is a normal subgroup of
$\mathrm{Aut}(\fg)$, we have
$f=h\circ\sigma\in I(\fg)$ if $h\in I(\fg)$. However
\begin{equation}\label{eq.18}
2\grad_{\tilde\fg}h(x+\lambda a)
=\grad_{\fg}h(x+\lambda a)
+\grad_{\fg}f(x-\lambda a)
\qquad\text{for any}\quad x\in\tilde\fg,
\end{equation}
because  $\sigma(a)=-a$, and $\sigma(x)=x$.
The five functions $h_1^{\lambda_1},h_2^{\lambda_2}$,
$f_1^{-\lambda_1},f_2^{-\lambda_2}$  and
$h_{a,b}$ commute pairwise on $\fg\simeq\fg^*$ with respect to
the standard (linear) Lie-Poisson bracket on $\fg$~\cite{MF78}.
This means that for any pair of functions $F_1,F_2$ from this set
we have
$$
\langle x,[\grad_{\fg} F_1(x), \grad_{\fg} F_2(x)]\rangle=0
\qquad\text{for all}\qquad
x\in\fg.
$$
Then by~\eqref{eq.18}
$$
\langle x,[\grad_{\tilde\fg} h_1^{\lambda_1}(x), \grad_{\tilde\fg}
h_2^{\lambda_2}(x)]\rangle=0
\qquad\text{for all}\qquad
x\in\tilde\fm\subset\fg.
$$
Now taking into account that
$(\grad_{\tilde\fg} h_j^{\lambda_j}(x))_{\tilde\fm}=
\grad_{\tilde\fm}\tilde H_j^{\lambda_j}(x)\in\tilde\fm(x)$,
$[x,\tilde\fm(x)]\subset\tilde\fm$ for $x\in\tilde\fm$
and $\tilde\fm\bot\tilde\fk$, we obtain that
$$
\langle x,[(\grad_{\tilde\fg}
h_1^{\lambda_1}(x))_{\tilde\fm}, (\grad_{\tilde\fg}
h_2^{\lambda_2}(x))_{\tilde\fm}]\rangle=0,
$$
i.e. $\{\tilde H_1^{\lambda_1},
\tilde H_2^{\lambda_2}\}^\mathrm{can}(x)=0$.
Similarly we can show that
$\{\tilde H^{\lambda},\tilde H_{a,b}\}^\mathrm{can}=0$.
\qed\end{proof}

As follows from the lemma above the set
$\mathcal{F}(\fg,\tilde\fm)$ is an involutive subset of
$(A^{\tilde K}_{\tilde\fm},\tilde\eta^\mathrm{can})$.
Let us define the numbers
$q(\tilde\fm)$, $p(\tilde\fm)$,
$r(\tilde\fm)$ and the subset
$R(\tilde\fm)\subset\tilde\fm$ similarly to the numbers
$q(\fm)$, $p(\fm)$, $r(\fm)$ and the subset
$R(\fm)\subset\fm$ but for the pair of algebras
$(\tilde\fg,\tilde\fk)$ (see~\eqref{eq.4}).

Let $\tilde\mu^\mathrm{can}: T(\tilde G/\tilde K)\to\tilde\fg^*$
be the standard moment map associated with the natural action of
$\tilde G$ on
$(T(\tilde G/\tilde K),\tilde\Omega)$). The subset
$\{\tilde h\circ\tilde\mu^\mathrm{can}, \, \tilde
h\in\mathcal{E}(\tilde\fg^*)\}\subset \mathcal{E}(T(\tilde
G/\tilde K))$
is a Poisson subalgebra of
$(\mathcal{E}(T(\tilde G/\tilde K)),\tilde\eta^\mathrm{can})$.
Since the Lie group
$\tilde G$ is compact, this subalgebra contains some
complete involutive subset of functions
$\mathcal{H}(\tilde\fg^*)$~\cite{Br86,My83}.
By definition of the moment map $\tilde\mu^\mathrm{can}$,
the union of the sets $\mathcal{H}(\tilde\fg^*)$ and
$\mathcal{F}(\fg,\tilde\fm)\subset A^{\tilde K}_{\tilde\fm}=A^{\tilde G}$
is an involutive subset of the Poisson algebra
$(\mathcal{E}(T(\tilde G/\tilde K)),\tilde\eta^\mathrm{can})$.

\begin{theorem}\label{th.10}
Suppose that the Zariski open subset $O^\mathrm{Kr}(\fm)\cap \tilde\fm$
of $\tilde\fm$ is nonempty.
Then the set $\mathcal{F}(\fg,\tilde\fm)$ is a complete involutive
subset of the algebra $(A^{\tilde K}_{\tilde\fm},\tilde\eta^\mathrm{can})$.
Moreover, the set $\mathcal{H}(\tilde\fg^*)\cup \mathcal{F}(\fg,\tilde\fm)$
is a complete involutive subset of the algebra
$(\mathcal{E}(T(\tilde G/\tilde K)),\tilde\eta^\mathrm{can})$
and the functions from this set
are integrals for the Hamiltonian flow with the Hamiltonian function
$\tilde H_{a,b}$ on $T(\tilde G/\tilde K)$.
\end{theorem}
\begin{proof}[Proof]
Fix some point $x\in O^\mathrm{Kr}(\fm) \cap R(\tilde\fm)\subset\tilde\fm$.
Put
\begin{equation*}
\begin{split}
L_x&=\{\grad_\fm f(x), f\in\mathcal{F}(\fg,\fm)\}
\subset\fm(x)\subset\fm,\\
\tilde L_x&=\{\grad_{\tilde\fm} \tilde f(x),
\tilde f\in \mathcal{F}(\fg,\tilde\fm)\}
\subset\tilde\fm(x)\subset\tilde\fm.
\end{split}
\end{equation*}
It is evident that $\tilde L_x=(L_x)_{\tilde\fm}$, where
$(\cdot)_{\tilde\fm}$ denotes the projection onto
$\tilde\fm$ along $\fm'$ in
$\fm=\tilde\fm\oplus\fm'$. Moreover, since
$x\in\tilde\fm$ and $\fk=\tilde\fk\oplus\fk'$, we have
$\ad x(\fk)=\ad x(\tilde\fk)\oplus\ad x(\fk')$, where
$\ad x(\tilde\fk)\subset\tilde\fm$ and
$\ad x(\fk')\subset\fm'$, and therefore
\begin{equation}\label{eq.19}
\fm(x)=\tilde\fm(x)\oplus(\fm(x)\cap\fm').
\end{equation}
By Lemma~\ref{le.9} the space
$\tilde L_x$ is an isotropic subspace of
$\tilde\fm(x)$ with respect to the form
$\tilde B_x: (y_1,y_2)\mapsto -\langle x,[y_1,y_2]\rangle$
on $\tilde\fm(x)$ associated with Poisson
bracket~\eqref{eq.17}. Therefore, in order to prove the
theorem it is sufficient to show that this subspace is
maximal isotropic.

To this end suppose that
$\langle x,[y,\tilde L_x]\rangle=0$ for some
$y\in\tilde\fm(x)$. Then $\langle [x,y],\tilde L_x\rangle=0$
by invariance of the form
$\langle\cdot ,\cdot \rangle$. Taking into account that
$[x,\tilde\fm(x)]\subset\tilde\fm$
by definition and
$\tilde\fm\bot\fm'$, we obtain that
$$
0=\langle [x,y],\tilde L_x\rangle
=\langle [x,y], L_x\rangle
=\langle x,[y, L_x]\rangle.
$$
But $y\in\fm(x)$ by~\eqref{eq.19}. Also by Theorem~\ref{th.7}
the space $L_x$ is a maximal isotropic in
$\fm(x)$ and, consequently, $y\in L_x$. Then
$y\in\tilde L_x$, because $y\in\tilde\fm$ and
$L_x\cap\tilde\fm\subset (L_x)_{\tilde\fm}=\tilde L_x$.
In other words, $\tilde L_x$ is a maximal isotropic subspace in
$\tilde\fm(x)$ with respect to the form
$\tilde B_x$ defined by~(\ref{eq.17}) on
$\tilde\fm(x)$.

Therefore the set
$\mathcal{F}(\fg,\tilde\fm)$ is a complete
involutive subset of the algebra
$A^{\tilde K}_{\tilde \fm}=A^{\tilde G}$
with respect to the canonical Poisson structure on
$(T(\tilde G/\tilde K),\tilde\Omega)$.
Then by~\cite[Lemma 3]{BJ01}
(see also~\cite[\S 2]{My01}, \cite{Pana03}) the set
$\mathcal{H}(\tilde\fg^*)\cup \mathcal{F}(\fg,\tilde\fm)$
is a complete involutive subset of the Poisson algebra
$(\mathcal{E}(T(\tilde G/\tilde K)),\tilde\eta^\mathrm{can})$
on the manifold $(T(\tilde G/\tilde K),\tilde\Omega)$.
By Lemma~\ref{le.9} the functions from the set $\mathcal{F}(\fg,\tilde\fm)$
are integrals of the Hamiltonian flow with the Hamiltonian function
$\tilde H_{a,b}$ on $T(\tilde G/\tilde K)$.
The functions from the set $\mathcal{H}(\tilde\fg^*)$
are also integrals for this flow because
$\tilde H_{a,b}\subset A^{\tilde G}$.
\qed
\end{proof}

\begin{remark}\label{re.11}
First of all note that the Hamiltonian flow of the function
$\tilde H_{a,a}$ (with $b=a$) coincides with the geodesic
flow determined by the Riemannian metric
$\langle\cdot ,\cdot\rangle$ on
$\tilde G/\tilde K$. The Hamiltonian function
$\tilde H_{a,b}$ on $T\tilde{\mathcal O}$
is the restriction of the Hamiltonian function
$H_{a,b}(x)=\frac12\langle x,\varphi_{a,b}(x)\rangle$,
$x\in\fm$ (considered as a $G$-invariant function on
$T{\mathcal O}$ from the set
$A^{G}=A^{K}_{\fm}$). Moreover, the function
$H_{a,b}$ is associated with the (pseudo)-Riemannian metric
$g_{a,b}$ on
${\mathcal O}$ (defined by the symmetric bilinear form
$(x,y)\mapsto\langle x,\varphi_{a,b}^{-1}(y)\rangle$,
$x,y\in\fm$) and $\tilde H_{a,b}$ with its restriction
$\tilde g_{a,b}=g_{a,b}|_{\tilde{\mathcal O}}$. Since
$\grad_\fm H_{a,b}(x)=\grad_{\tilde\fm}\tilde
H_{a,b}(x)=\varphi_{a,b}(x)\in\tilde\fm$
for each
$x\in\tilde\fm$, it is easy to verify that the Hamiltonian
vector field $X_{H_{a,b}}$ of $H_{a,b}$ on
$T{\mathcal O}$ is tangent to the submanifold
$T\tilde{\mathcal O}\subset T{\mathcal O}$
(using, for instance, the exact expression for
$X_{H_{a,b}}$ on $T{\mathcal O}$~\cite[(3.6)]{MP04}).
As follows from Proposition~\ref{pr.4},
$(T\tilde{\mathcal O},\tilde\Omega)$
is a symplectic submanifold of
$(T{\mathcal O},\Omega)$. Therefore the Hamiltonian vector
field $X_{\tilde H_{a,b}}$ of $\tilde H_{a,b}$ on
$T\tilde{\mathcal O}$ coincides with the restriction
$X_{H_{a,b}}|_{T\tilde{\mathcal O}}$
and, consequently, the suborbit
$(\tilde{\mathcal O},\tilde g_{a,b})$
is a totally geodesic submanifold of the orbit
$({\mathcal O},g_{a,b})$. Taking into account that the
geodesic flow of the metric
$g_{a,b}$ on $T{\mathcal O}$ is integrable
(see~\cite[Theorem 3.10]{MP04} or~\cite{BJ04}),
Theorem~\ref{th.10} above is an illustration
to the phenomenon of a complete integrability of the
geodesic flow on a totally geodesic submanifold of a
manifold with the integrable geodesic flow.
\end{remark}

\begin{remark}\label{re.12}
Let $\fg_1$ be some $\sigma$-invariant subalgebra of
$\fg$ containing its maximal semisimple ideal
$\fg_s$. Assume that the connected subgroup
$G_1\subset G$ with the Lie algebra
$\fg_1$ is closed. Denote by
$a_1\in\fg_1$ the orthogonal projection of the element
$a\in\fg$. Put $\fk_1=\fg_1^{a_1}$
($\fk=\fg^a$). Since the algebra
$\fk=\fg^a$ contains the center $\fz(\fg)$ of
$\fg$, the space $\fm=\fk^\bot$ is a subspace of the ideal
$\fg_s$ and the space
$\fm_1=\fk_1^\bot\subset\fg_1$ coincide with
$\fm$. Consequently, the bi-Poisson structures
$\eta^t(\fg,a)$ and $\eta^t(\fg_1,a_1)$
are either Kronecker or non-Kronecker simultaneously at a point
$x\in R(\fm)\cap\tilde\fm$ if
$R(\fm)\cap\tilde\fm\ne\emptyset$.
\end{remark}

\subsection{Conditional reduction theorem}
\label{ss.2.4}

We will use the notation introduced in the previous
subsection. Suppose that the set
$R(\fm)\cap\tilde\fm$ is nonempty and that
$x_0$ is a common element of this set and the set
$R(\tilde\fm)$. Let $\fg_0$ be the centralizer of
$\fk^{x_0}$ in $\fg$. Put $\fk_0=\fk\cap\fg_0$ and
$\fm_0=\fm\cap\fg_0$. Since $\sigma(x_0)=x_0$ and
$\sigma(\fk)=\fk$, we have
$\sigma(\fk^{x_0})=\fk^{x_0}$. Therefore the spaces
$\fg_0$, $\fk_0$ and $\fm_0$ are
$\sigma$-invariant, i.e.
$$
\fg_0=\fk_0\oplus\fm_0,
\quad
\tilde\fg_0=\tilde\fk_0\oplus\tilde\fm_0,
\quad\text{where}\
\tilde\fg_0=\fg_0\cap\tilde\fg,\
\tilde\fk_0=\fk_0\cap\tilde\fk,\
\tilde\fm_0=\fm_0\cap\tilde\fm.
$$
By Proposition~\ref{pr.5}
$x_0\in R(\fm_0)$, $a\in\fg_0$ and $\fk_0=\fg_0^a$
($\fk=\fg^a$). Now we will prove the following conditional
reduction theorem which reduces the pair
$(\fg,\fk=\fg^a)$ to the pair
$(\fg_0,\fk_0=\fg_0^{a})$ if the condition
$R(\fm)\cap\tilde\fm\ne\emptyset$ holds.
To this end, define the sets $O^\mathrm{Kr}(\fm_0)$,
$Q_a(\fm_0)$ and $M_a(\fm_0)$
similarly to the sets $O^\mathrm{Kr}(\fm)$,
$Q_a(\fm)$ and $M_a(\fm)$
but for the pair of algebras
$(\fg_0,\fk_0)$ (see~(\ref{eq.14}), (\ref{eq.11}) and~(\ref{eq.13})).

\begin{theorem}\label{th.13}
Suppose that
$R(\fm)\cap\tilde\fm\ne\emptyset$ and choose any point
$x_0\in R(\fm)\cap R(\tilde\fm)$. The set
$O^\mathrm{Kr}(\fm)\cap \tilde\fm$
is nonempty if so is the set
$O^\mathrm{Kr}(\fm_0)\cap \tilde\fm_0$. In particular, if
$O^\mathrm{Kr}(\fm_0)\cap \tilde\fm_0\ne\emptyset$
then the set
$\mathcal{F}(\fg,\tilde\fm)$ is a complete involutive subset
of the algebra
$(A^{\tilde K}_{\tilde\fm},\tilde\eta^\mathrm{can})$.
\end{theorem}
\begin{proof}[Proof]
Assume that $O^\mathrm{Kr}(\fm_0)\cap \tilde\fm_0\ne\emptyset$.
By definition, $x_0\in R(\fm)\cap R(\tilde\fm)$ and
$x_0\in\tilde\fm_0$ by the construction above. Therefore the set
$R(\fm)\cap R(\tilde\fm)\cap O^\mathrm{Kr}(\fm_0)\cap \tilde\fm_0$
is a nonempty Zariski open subset of
$\tilde\fm_0$. Thus the set
$R(\fm)\cap R(\tilde\fm)\cap O^\mathrm{Kr}(\fm_0)\cap R(\tilde\fm_0)$
is also nonempty. Fix some element
$x$ belonging to this set. To prove the theorem it is
sufficient to show that
$x\in O^\mathrm{Kr}(\fm)$.

Since $x\in R(\fm)\cap \fm_0$, by Proposition~\ref{pr.5}
$(1)$ we have that $\fm(x)=\fm_0(x)$.
Then by relations~(\ref{eq.10}) the two bi-Poisson structures
$\eta(\fg,a)$ and $\eta(\fg_0,a)$ determine on the space
$\fm(x)=\fm_0(x)$ the same families of skew-symmetric
operators. In other words, at the point
$x\in R(\fm)\cap R(\fm_0)$ these bi-Poisson structures are
Kronecker simultaneously. Thus
$x\in O^\mathrm{Kr}(\fm)\cap \tilde\fm$. \qed
\end{proof}

As follows from Theorem~\ref{th.10} we have to establish
when the set $O^\mathrm{Kr}(\fm)\cap\tilde\fm$ is nonempty.
Since the subsets $Q_a(\fm)\subset\fm$ and $M_a(\fm)\subset\fm$
are Zariski open (see~(\ref{eq.12}) and Lemma~\ref{le.8}),
we obtain by~(\ref{eq.14}) that
\begin{equation}\label{eq.20}
O^\mathrm{Kr}(\fm)\cap\tilde\fm\ne\emptyset
\qquad\iff\qquad Q_a(\fm)\cap\tilde\fm\ne\emptyset
\quad\text{and}\quad M_a(\fm)\cap\tilde\fm\ne\emptyset.
\end{equation}

By Proposition~\ref{pr.5}
$(4)$ for any element
$x\in R(\fm)\cap\fm_0$ its centralizer
$\fk_0^{x_0}$ in $\fk_0$ is the center
$\fz(\fg_0)$ of the Lie algebra $\fg_0$.
As follows from Theorem~\ref{th.13} and
relations~(\ref{eq.20}) to prove a completeness of the set
$\mathcal{F}(\fg,\tilde\fm)$ it is sufficient to show that
$Q_a(\fm_0)\cap\tilde\fm_0\ne\emptyset$ and
$M_a(\fm_0)\cap\tilde\fm_0\ne\emptyset$.
Therefore we will investigate the two open subsets
$Q_a(\fm)\cap\tilde\fm$ and
$M_a(\fm)\cap\tilde\fm$ in more detail in the following two
subsections only in the case when for each element
$x\in R(\fm)$ its centralizer
$\fk^x=\fg^a\cap\fg^x$ is the center $\fz(\fg)$ of
$\fg$.

\subsection{Necessary and sufficient conditions for the
set $Q_a(\fm)\cap\tilde\fm$ to be nonempty}
\label{ss.2.5}

We will use the notation introduced in subsection~\ref{ss.2.3}.
\begin{proposition}\label{pr.14}
Suppose that
$R(\fm)\cap\tilde\fm\ne\emptyset$
and choose an arbitrary point
$x_0\in R(\fm)\cap R(\tilde\fm)$.
Let $\fk^{x_0}=\fz(\fg)$.
The set $Q_a(\fm)\cap\tilde\fm$ is nonempty
if and only if
\begin{equation}\label{eq.21}
\text{the subspace } \fk'=(1-\sigma)\fk
\text{ contains regular elements of the Lie algebra }
\fk.
\end{equation}
\end{proposition}
\begin{proof}[Proof]
Let $V$ be some vector subspace of the space $\fm$
for which $V\cap R(\fm)\ne\emptyset$. Put
$$
m_a(V)=\min_{x\in V\cap R(\fm)}
\dim \Bigl(\fm(x)\cap \bigl(\ad_a^{-1}\ad x(\fk)\bigr)\Bigr),
$$
where, recall, $\ad_a^{-1}\eqdef (\ad a|_\fm)^{-1}$.
By~\eqref{eq.12},
$m_a(\fm)=r(\fm)$ and clearly
$m_a(V)\geqslant r(\fm)$. We have to show that
$m_a(\tilde\fm)=r(\fm)$ if and only if
condition~(\ref{eq.21}) holds.

We will prove this proposition using the moment map theory
and the method proposed in~\cite{Pana03}. Since
$\fk^{x_0}=\fz(\fg)$, then by the dimension arguments
$\fk^{x}=\fz(\fg)$ for all $x\in R(\fm)$.
For our aim it is convenient to use the moment map
constructed in our previous paper~\cite{MP04}. To this end,
here we briefly describe main properties of this moment
map~\cite[Remark 3.2]{MP04}.

Consider on the vector space $\fm$ the non-degenerate
bilinear form
$$
\beta(y_1,y_2)=\langle
y_1,\ad_a^{-1}(y_2)\rangle,\qquad y_1,y_2\in\fm.
$$
Since the endomorphism
$\ad a|_\fm:\fm\to\fm$ is skew-symmetric (with respect to the
form $\langle\cdot,\cdot\rangle$), the form $\beta$
is also skew-symmetric. Identifying the tangent space
$T_x\fm$ with $\fm$ for each $x\in\fm$, we can consider
$\beta$ as a symplectic form on $\fm$. Since the
$\Ad$-action of $K$ on $\fm$ preserves the form
$\beta$, this action of $K$ is Hamiltonian with the
$K$-equivariant moment map
$$\mu^\beta:\fm\to\fk^*,
\quad
\mu^\beta(x)(\zeta)=-\frac12\langle \ad_a^{-1} (x),
[\zeta,x]\rangle, \quad \forall \zeta\in\fk
$$
(see~\cite[Remark 3.2]{MP04}). In
particular, for each $\zeta\in\fk$ the vector field
$\zeta_X(x)=[\zeta,x]\in T_x\fm$
is the Hamiltonian vector field of the function
$f_\zeta(x)=\mu^\beta(x)(\zeta)$ on the manifold
$X=\fm$.

Let $x\in\fm$, $W_x\subset T_x\fm$
be the tangent space to the
$K$-orbit in $\fm$ and let
$W_x^{\beta\bot}$ be the (skew)orthogonal complement to
$W_x$ in $T_x\fm$ with respect to the form
$\beta$. It is easy to see that $W_x=\ad x(\fk)$
and $W_x^{\beta\bot}=\ad a (\fm(x))$, i.e.
$$
\dim\bigl(W_x\cap W_x^{\beta\bot}\bigr)
=\dim\Bigl(\ad a\bigl(\fm(x)\bigr)\cap \ad x(\fk)\Bigr)
=\dim \Bigl(\fm(x)\cap \bigl(\ad_a^{-1}\ad x(\fk)\bigr)\Bigr).
$$
However, by the $K$-equivariance of the moment map $\mu^\beta$,
$\zeta_X(x)\in W_x\cap W^{\beta\bot}_x$ if and only if
$\ad^*\zeta(\alpha)=0$, where
$\alpha=\mu^\beta(x)$~\cite{GS83}.
In other words, the dimension
$\dim\bigl(W_x\cap W_x^{\beta\bot}\bigl)+\dim\fz(\fg)$
equals the codimension of the orbit
$\Ad^*(K)\cdot\alpha$ in $\fk^*$ for each
$x\in R(\fm)$ ($\zeta_X(x)=0$ if and only if
$\zeta\in\fk^x=\fz(\fg)$).

Identifying the space $\fk$ with its dual $\fk^*$
by means the form $\langle\cdot ,\cdot \rangle$, we obtain that
$$
\mu^\beta:\fm\to\fk, \qquad
\mu^\beta(x)=\frac12\,[\ad_a^{-1} (x),x]_\fk
$$
and
$\dim\bigl(W_x\cap W_x^{\beta\bot}\bigr)=
\dim \fk^\alpha-\dim\fz(\fg)$ if $x\in R(\fm)$. Thus
$$
m_a(V)=\min_{x\in V\cap R(\fm)} \dim \fk^{\alpha(x)}-\dim\fz(\fg)
=\min_{x\in V}\dim \fk^{\alpha(x)}-\dim\fz(\fg),
$$
where $\alpha(x)=\mu^\beta(x)\in\fk$ and $V\cap R(\fm)\ne\emptyset$.

Since $\fz(\fg)\subset\fk$, the space $\fm$
is contained in the maximal semi-simple ideal $\fg_s$
of $\fg$. Remark also that
$\mu^\beta(\tilde\fm)\subset[\fg,\fg]_\fk\subset\fg_s$ and
by~(\ref{eq.16})
$\mu^\beta(\tilde\fm)\subset [\fm',\tilde\fm]_\fk\subset(\fg')_\fk=\fk'$.

To determine the number
$\min_{x\in \tilde\fm} \dim \fk^{\alpha(x)}$
we will show that the image
$\mu^{\beta}(\tilde\fm)$ contains an open subset in the
space $\fk'\cap\fg_s$. It is easy to calculate that, for any
tangent vector $y\in\fm=T_{x_0}\fm$,
$$
D_{x_0}(y)\eqdef\mu^{\beta}_*(x_0)(y)=\frac12[\ad^{-1}_a (x_0),y]_{\fk}
+\frac12[\ad^{-1}_a (y),x_0]_{\fk}.
$$
Taking into account relations~\eqref{eq.15} and~\eqref{eq.16}
and the inclusion $x_0\in\tilde\fm\subset\tilde\fg$,
we obtain that
$$
D_{x_0}(\fm')\subset \left([\ad^{-1}_a(\tilde\fm),\fm']
+[\ad^{-1}_a(\fm'),\tilde\fm]\right)_{\fk}
\subset(\tilde\fg)_{\fk}\subset\tilde\fk.
$$
and, similarly,
$D_{x_0}(\tilde\fm)\subset(\fg')_{\fk}\subset\fk'$.
In other words, $D_{x_0}(\fm)=D_{x_0}(\fm')\oplus
D_{x_0}(\tilde\fm)$.

The image $\mu^{\beta}_*(T_{x_0}\fm)\subset\fk$
of the tangent map of the moment map
$\mu^{\beta}$ at $x_0$ coincides with the annihilator in
$\fk^*\simeq\fk$ of the Lie algebra
$\fk^{x_0}$ of the isotropy group
$\{k\in K_0: \Ad k(x_0)=x_0\}$ of
$x_0\in\fm$~\cite{GS83}. Since this annihilator coincides
with the orthogonal complement of the center
$\fz(\fg)$ in $\fk$, we obtain that
$D_{x_0}(\fm)=\fk\cap\fg_s$. But
$\fk\cap\fg_s=(\fk'\cap\fg_s)\oplus(\tilde\fk\cap\fg_s)$
because $\sigma(\fg_s)=\fg_s$. Therefore
$\mu^{\beta}_{*}(x_0)(\tilde\fm)=\fk'\cap\fg_s$
and, consequently, the set
$\mu^{\beta}(\tilde\fm)$ contains some open subset in
$\fk'\cap\fg_s$. Also $\fk^\alpha=\fk^{\alpha+z}$ for any
$\alpha\in\fk'$, $z\in\fz(\fg)$. Therefore
\begin{equation}\label{eq.22}
m_a(\tilde\fm)
=\min_{\alpha\in \fk'\cap\fg_s} \dim \fk^\alpha-\dim\fz(\fg)
=\min_{\alpha\in \fk'} \dim \fk^\alpha-\dim\fz(\fg).
\end{equation}
Since $\fk^{x_0}=\fz(\fg)$ is a commutative algebra, each element
of the set $R(\fm)$ is a regular element of $\fg$, i.e.
$q(\fm)=\rank\fg$ and
$r(\fm)=\rank\fg-\dim\fz(\fg)$ (see~(\ref{eq.7})). Also
$\rank\fk=\rank\fg$ because
$\fk=\fg^a$. Hence by~(\ref{eq.22})
$m_a(\tilde\fm)=r(\fm)$ if and only if the algebra
$\fk'$ contains regular elements of the Lie algebra
$\fk$. \qed
\end{proof}

\subsection{Sufficient conditions for the
set $M_a(\fm)\cap\tilde\fm$ to be nonempty}
\label{ss.2.6}

In this subsection we will give a sufficient condition for
an open subset $M_a(\fm)\cap\tilde\fm$ of $\tilde\fm$
to be nonempty. This condition is motivated by some trick
used in the proof of Theorem~3.4 in~\cite{BJ04}. Here we
will use the notation introduced in subsection~\ref{ss.2.3}.

Let $\ft$ be a Cartan subalgebra of
$\fg$ containing the element
$a$. The complex space $\ft^\bbC$
is a Cartan subalgebra of the reductive complex Lie algebra
$\fg^\bbC$. Since
$\fk=\fg^a$ and the commutative Lie algebra
$\ft$ contains $a$,
$\ft^\bbC$ is also a Cartan subalgebra of
$\fk^\bbC$. Let $\Delta$ be the root system of
$\fg^\bbC$ with respect to $\ft^\bbC$. Denote by
$\Delta_\fk$ the set of roots in
$\Delta$ which vanish identically on
$a$. This is the root system of $(\fk^\bbC,\ft^\bbC)$. Put
$\Delta_\fm=\Delta\setminus\Delta_\fk$.
Then we have the direct decompositions
$$
\fg^\bbC=\ft^\bbC\oplus\sum_{\alpha\in\Delta}
\fg^\alpha,
\qquad
\fk^\bbC=\ft^\bbC\oplus
\sum_{\alpha\in\Delta_\fk}\fg^\alpha,
\qquad
\fm^\bbC=\sum_{\alpha\in\Delta_\fm}\fg^\alpha,
$$
where $\fg^\alpha$ is the root space corresponding to the root
$\alpha$. We can choose a Weyl basis $\{E_\alpha,\alpha\in\Delta\}$
of $\fg^\bbC$ mod $\ft^\bbC$ such that
$E_\alpha\in\fg^\alpha$ for each $\alpha\in\Delta$ and
the compact form $\fg$ of $\fg^\bbC$ is spanned by
$\ft$ and the vectors $X_\alpha=(E_\alpha-E_{-\alpha})$,
$Y_\alpha=i(E_\alpha+E_{-\alpha})$, $\alpha\in\Delta$.

\begin{proposition}\label{pr.15}
Suppose that
$R(\fm)\cap\tilde\fm\ne\emptyset$
and $\fk^{x}=\fz(\fg)$ for all $x\in R(\fm)$.
Assume that there exists a system of simple roots
$\pi$ of $\Delta$ such that $\pi\subset\Delta_\fm$
and some vector
$$
x_\pi=\sum_{\alpha\in\pi} c_\alpha E_{-\alpha}
+\sum_{\beta\in\Delta_\fm^+}d_\beta E_{\beta},
\qquad
c_\alpha\in\bbC\setminus\{0\}, \alpha\in\pi,
\ d_\beta\in\bbC, \beta\in\Delta_\fm^+,
$$
where $\Delta^+$ is the system of positive roots of
$\Delta$ determined by $\pi$
and $\Delta_\fm^+=\Delta_\fm\cap \Delta^+$, belongs to the space
$\tilde\fm^\bbC\subset\fm^\bbC$. Then the set
$M_a(\fm)\cap\tilde\fm$ is nonempty.
\end{proposition}
\begin{proof}[Proof]
Since for $x\in R(\fm)$ the algebra $\fk^{x}=\fz(\fg)$
is commutative, the centralizer $\fg^x$ is a Cartan subalgebra of $\fg$
(see~(\ref{eq.5})), i.e. $q(\fm)=\rank\fg$.

The element $e_-=\sum_{\alpha\in\pi} c_\alpha E_{-\alpha}$
is a principal nilpotent element of the reductive complex Lie algebra
$\fg^\bbC$~\cite{Ko63}. Since $a\in\ft^\bbC$, for each
$\lambda\in\bbC$ the vector
$x_\pi+\lambda a$ is an element of the affine subspace
$e_- +\mathfrak{b}$, where
$\mathfrak{b}=\ft^\bbC\oplus\sum_{\beta\in\Delta^+} \fg^\beta$
is a Borel subalgebra of
$\fg^\bbC$. Then by~\cite[Lemma~10]{Ko63}
$\dim_\bbC(\fg^\bbC)^{x_\pi+\lambda a}=\rank\fg^\bbC$. But
$q(\fm)=\rank\fg^\bbC$. Thus
$x_\pi\in M_a(\fm^\bbC)\cap\tilde\fm^\bbC$
(see the proof of Lemma~\ref{le.8}). Since
$M_a(\fm^\bbC)\cap\tilde\fm^\bbC$
is a nonempty Zariski open subset of
$\tilde\fm^\bbC$, its intersection with the real form
$\tilde\fm\subset\tilde\fm^\bbC$
is nonempty. Taking into account that
$M_a(\fm^\bbC)\cap\tilde\fm=M_a(\fm)\cap\tilde\fm$
we complete the proof. \qed
\end{proof}

\subsection{ Integrable geodesic flows on
$SO(n)/\bigl(SO(n_1)\times\cdots\times SO(n_p)\bigr)$}
\label{ss.2.7}

In this subsection we show that the conditions of
Theorem~\ref{th.10} hold for the homogeneous space
$SO(n)/\bigl(SO(n_1)\times\cdots\times SO(n_p)\bigr)$,
$n_1+\cdots+n_p=n$. The case $p=2$
is not interesting for us from the point of view of
integrability because in this case the considered homogeneous space
$\tilde G/\tilde K$ is a symmetric space and, consequently,
all $\tilde G$-invariant Hamiltonian flows on
$T(\tilde G/\tilde K)$ are integrable~\cite{Mi82,My86}.

Consider the symmetric space
$G/\tilde G=U(n)/SO(n)$, where $n\geqslant4$, with the involution
$\sigma$ on the Lie algebra of skew-hermitian
matrices $\fg=\mathfrak{u}(n)$ determined by the complex conjugation. Then the
Lie algebra $\tilde\fg=(1+\sigma)\fg$ is the Lie algebra
$\mathfrak{so}(n)$ of all real skew-symmetric
$n\times n$ matrices. The space $\fg'=(1-\sigma)\fg$
coincides with the set $i\,\mathrm{sym}(n)$, where
$\mathrm{sym}(n)$ is the space of all real symmetric
$n\times n$ matrices.

Fix some element $a\in\fg'$,
$a=\mathrm{diag}(i\lambda_1,...,i\lambda_1,
i\lambda_2,...,i\lambda_2,\ldots,i\lambda_p,...,i\lambda_p)$,
where all real numbers
$\lambda_1,...,\lambda_p$ are pairwise different and the
multiplicity of each $i\lambda_j$ is equal to $n_j\geqslant1$,
$n_1+\cdots+n_p=n$. Without loss of generality (to simplify
calculations) we may assume that
$1\leqslant n_1\leqslant n_2\leqslant\cdots\leqslant n_p<n$.

It is clear that the Lie algebra
$\fk=\fg^a$ is the Lie algebra
$\mathfrak{u}(n_1)\oplus\cdots\oplus\mathfrak{u}(n_p)$
(with the standard block-diagonal embedding) and
$\tilde\fk$ is the real part of this Lie algebra, i.e.
$\tilde\fk$ coincides with
$\mathfrak{so}(n_1)\oplus\cdots\oplus\mathfrak{so}(n_p)$
($\mathfrak{so}(1)=0$). In this case,
$\tilde G/\tilde K=SO(n)/\bigl(SO(n_1)\times\cdots\times
SO(n_p)\bigr)$.

Putting $\langle X,Y \rangle=-\mathrm{Tr}XY$
(using the trace-form) we define an invariant scalar product on
$\fg$. To describe the space
$\fm=\fk^\bot$ consider any matrix
$X\in\fg$ as a block-matrix consisting of rectangle elements
$X^{k,l}$, which are rectangle complex
$n_k\times n_l$ matrices,
$1\leqslant k,l\leqslant p$. It is clear that
$(\overline{X^{k,l}})^t=-X^{l,k}$
and therefore any element of
$\mathfrak{u}(n)$ is defined by its blocks $X^{k,l}$ with
$k\leqslant l$. As a space, the Lie algebra
$\fg=\mathfrak{u}(n)$ is a direct sum of its block-type
subspaces $V^{k,l}$, $1\leqslant k\leqslant l\leqslant p$.
In this notation the Lie subalgebra
$\fk$ is the direct sum $\sum_{j=1}^p V^{j,j}$ and
$\fm=\sum_{1\leqslant k<l\leqslant p} V^{k,l}$.
We will denote the corresponding to
$X^{k,l}$ element of the space $V^{k,l}$ by
$\varphi(X^{k,l})$. Each subspace $V^{k,l}$ is an
$\fk$-module, i.e. $[V^{k,l},\fk]\subset V^{k,l}$.

First of all consider the simplest case when
$p=2$. In this case
$G/K=U(n_1+n_2)/\bigl(U(n_1)\times U(n_2)\bigr)$
is a Hermitian symmetric space. Therefore there exists a
Cartan subspace $\mathfrak{a}$ in
$V^{1,2}=\fm$ (a maximal commutative subspace in
$V^{1,2}$) consisting of real matrices (belonging to
$\mathfrak{so}(n)$)~\cite[Ch.X, sec.2.1]{He78}. This
$n_1$-dimensional Cartan subspace
$\mathfrak{a}$ can be described by the ``diagonal" matrices
$X^{1,2}$ in which all entries vanish with the exception of
$n_1$ entries $X^{1,2}_{j,j}$,
$j=1,..,n_1\leqslant n_2$, which are arbitrary real numbers.
Then the centralizer
$\fk^{x_0}$ of a regular element (of the Cartan subspace)
$x_0=\varphi(X^{1,2}_*)\in\tilde\fm\subset V^{1,2}$ in
$\mathfrak{u}(n_1)\oplus\mathfrak{u}(n_2)=V^{1,1}\oplus V^{2,2}$
is a direct sum $\ft_*\oplus\mathfrak{c}_*$, where
$\mathfrak{c}_*\simeq\mathfrak{u}(n_2-n_1)$ and
$\ft_*$ is a commutative algebra of dimension
$n_1$ consisting of diagonal matrices
$\mathrm{diag}(ix_1,..,ix_{n_1},ix_1,..,ix_{n_1},0,..,0)$,
$x_j\in\bbR$. Remark that the maximal semisimple ideal of
$\mathfrak{c}_*$ coincides with the maximal semisimple ideal of the
centralizer of $\ft_*$ in
$\mathfrak{u}(n_1)\oplus\mathfrak{u}(n_2)$.

Suppose now that
$p\geqslant3$. We will attempt to describe the centralizer
$\fk^{x_0}$ of some generic element
$x_0\in\tilde\fm$ constructing this element. To simplify our
calculations remark that each space
$V^{k,k}\oplus V^{l,l}\oplus V^{k,l}$,
$k<l$ is a Lie subalgebra of
$\mathfrak{u}(n)$ isomorphic to $\mathfrak{u}(n_k+n_l)$,
$V^{k,k}\oplus V^{l,l}\simeq \mathfrak{u}(n_k)\oplus
\mathfrak{u}(n_l)$
and $[V^{k,l},V^{k,l}]\subset V^{k,k}\oplus V^{l,l}$. But
$U(n_k+n_l)/\bigl(U(n_k)\times U(n_l)\bigr)$
is a Hermitian symmetric space and therefore we can use our
calculations for the case
$p=2$. Since each subspace $V^{k,l}$ is a
$\fk$-module, we will construct the element
$x_0$ selecting step by step its
$V^{k,l}$-entries. For our aim it is enough to consider the
submodule
$\sum_{j=1}^{p-1} V^{j,j+1}\oplus\sum_{j=1}^{p-2} V^{j,p}$
of $\fm$. Choosing in each $\fk$-module
$ V^{j,j+1}$ of the first component the ``diagonal" element
$\varphi(X^{j,j+1}_*)$ as above, we obtain that their common
isotropy algebra is the direct sum
$\ft_*\oplus\mathfrak{c}_*$, where
$\mathfrak{c}_*\simeq\mathfrak{u}(n_p-n_{p-1})$ and
$\ft_*$ is of commutative algebra of dimension
$n_{p-1}$, consisting of the diagonal matrices
$\mathrm{diag}(ix_1,..,ix_{n_1},ix_1,..,ix_{n_2},..,ix_1,..
,ix_{n_{p-1}},ix_1,..,ix_{n_{p-1}},0,..,0)$,
$x_j\in\bbR$, $1\leqslant j\leqslant n_{p-1}$.
Remark that the maximal semisimple ideal of
$\mathfrak{c}_*$ coincides with the maximal semisimple ideal of the
centralizer of $\ft_*$ in $\fk$. Now we consider
$V=\sum_{j=1}^{p-2} V^{j,p}$ as a
$\ft_*\oplus\mathfrak{c}_*$-module (not as
a $\fk$-module) of complex dimension
$N_2\times n_p$, $N_2=n_1+..+n_{p-2}$. Then
$V$ is the direct sum of $\ft_*\oplus\mathfrak{c}_*$-modules
$V^{(1)}\oplus V^{(2)}$, $V^{(1)}\bot V^{(2)}$, where
$V^{(1)}$ (of complex dimension
$N_2\times n_{p-1}$) is a trivial
$\mathfrak{c}_*$-module. Therefore the isotropy subalgebra
of a real generic point from $V^{(1)}$
is the algebra $\ft_{**}\oplus\mathfrak{c}_*$, where
$\ft_{**}\subset\ft_*$
is the one-dimensional subalgebra
consisting of the elements of $\ft_*$ such that
$x_1=x_2=..=x_{n_{p-1}}=\lambda$.
Considering the module
$V^{(2)}$ as the space $M_{N_2,n_p-n_{p-1}}$ of complex
$N_2\times (n_p-n_{p-1})$ matrices,
the $\ad$-representation of $\ft_{**}\oplus\mathfrak{c}_*$ in
$V^{(2)}$ is described as follows:
$\ad(\lambda,A)(B) =i\lambda B-BA$, where
$(\lambda,A)\in \bbR\oplus\mathfrak{u}(n_p-n_{p-1})\simeq
\ft_{**}\oplus\mathfrak{c}_*$ and
$B\in M_{N_2,n_p-n_{p-1}}\simeq V^{(2)}$.

If the number of rows
$N_2=n_1+..+n_{p-2}$ in
$B$ is not less than the number of columns
$(n_p-n_{p-1})$, then for any real matrix
$B\in M_{N_2,n_p-n_{p-1}}$ of the maximal rank
we have $i\lambda B-BA=0$ if and only if
$A=i\lambda$ (is a scalar matrix). Therefore in this case
$\fk^{x_0}=\fz(\fg)$ for some (real) matrix
$x_0\in\tilde\fm$.

If $N_2<n_p-n_{p-1}$, then for the real matrix
$B$ in which all entries vanish with the exception of
$N_2$ entries $B_{j,j}=1$, $j=1,..,N_2$, we obtain that
$i\lambda B-BA=0$ if and only if
$A$ is an element of a Lie algebra isomorphic to
$\mathfrak{u}(1)\oplus\mathfrak{u}(n_p-n_{p-1}-N_2)$,
i.e. for some (real) matrix
$x_0\in\tilde\fm$ we have
$\dim\fk^{x_0}=1+\dim\mathfrak{u}(n_p-n_{p-1}-N_2)$.

Suppose now that $p\geqslant3$ and
$n_1+..+n_{p-2}\geqslant n_p-n_{p-1}$.
We showed above that for some element
$x_0\in\tilde\fm$ its centralizer
$\fk^{x_0}$ is the one-dimensional center
$\fz(\fg)$ of $\fg=\mathfrak{u}(n)$. Since
$\fz(\fg)\subset \fk^{x}$ for each $x\in\fm$,
we have that $\dim\fk^{x_0}=p(\fm)$ and
$\dim\tilde\fk^{x_0}=p(\tilde\fm)$.
Clearly these two properties hold for any points
$x\in\tilde\fm$ which are sufficiently closed to $x_0$.
Therefore we can suppose without
losing generality that
$x_0\in R(\tilde\fm)$, i.e.
$\dim\tilde\fg^{x_0}=q(\tilde\fm)$. But
$q(\tilde\fm)=\rank(\tilde\fg)$ because the algebra
$\tilde\fk^{x_0}\subset\fz(\fg)$
is commutative (see~(\ref{eq.5})). Taking into account that
any regular element of the Lie algebra
$\tilde\fg=\mathfrak{so}(n)$ is regular in
$\fg=\mathfrak{u}(n)$ we obtain that $x_0\in R(\fm)$. Thus
$x_0\in R(\fm)\cap R(\tilde\fm)$ and
$\fk^{x_0}=\fz(\fg)$.

Now we can use Proposition~\ref{pr.14} to prove that
$Q_a(\fm)\cap\tilde\fm\ne\emptyset$. Indeed, since the space
$\fk'=(1-\sigma)\fk=i\mathrm{sym}(n_1)\oplus\cdots\oplus
i\mathrm{sym}(n_p)$
contains regular elements of the Lie algebra
$\fk=\mathfrak{u}(n_1)\oplus..\oplus \mathfrak{u}(n_p)$,
by this proposition $Q_a(\fm)\cap\tilde\fm\ne\emptyset$.

Our element $a\in\mathfrak{u}(n)$ is a diagonal matrix
$\mathrm{diag}(a_1,..,a_n)$,
$a_j\in i\bbR$ and is contained in the Cartan subalgebra $\ft$ of
$\mathfrak{u}(n)$ consisting of matrices
$h=\mathrm{diag}(h_1,..,h_n)$, $h_j\in i\bbR$. The set
$\Delta=\{\varepsilon_{jk}=\varepsilon_j-\varepsilon_k,\
j\ne k; j,k=1,..,n\}$,
$\varepsilon_j(h)=h_j$, is the standard system of roots of
$\fg^\bbC$ with respect to $\ft^\bbC$. The algebra
$\ft$ and vectors
$X_{jk}=E_{\varepsilon_{jk}}-E_{\varepsilon_{kj}}$,
$Y_{jk}=iE_{\varepsilon_{jk}}+iE_{\varepsilon_{kj}}$ span
$\fg$ (by $E_{\varepsilon_{jk}}$
we denote the matrix with only nonzero element
$1$ occupies the position $jk$). To prove that
$M_a(\fm)\cap\tilde\fm\ne\emptyset$
we will use the following observation (see~\cite{BJ03}): since
$n_p\leqslant n_1+\cdots+n_{p-1}$,
there exists a permutation
$p$ of $\{1,2,..,n\}$ such that $a_{p(j)}\ne a_{p(j+1)}$,
$j=1,..,n-1$. In other words,
$(\varepsilon_j-\varepsilon_{j+1})(\hat a)\ne0$ for each
$j=1,..,n-1$, where
$\hat a=\mathrm{diag}(a_{p(1)},a_{p(2)},..,a_{p(n)})$.
As the Weyl group of $(\fg,\ft)$ is the permutation group of
$n$ elements, there exists the system of simple roots
$\pi$ of $\Delta$ such that $\alpha(a)\ne0$ for each
$\alpha\in\pi$. Since the element
$\sum_{\alpha\in\pi}(E_\alpha-E_{-\alpha})$
(real skew-symmetric matrix) belongs to the space
$\tilde\fm$, then by Proposition~\ref{pr.15},
$M_a(\fm)\cap\tilde\fm\ne\emptyset$ and by~(\ref{eq.20})
$O^\mathrm{Kr}(\fm)\cap\tilde\fm\ne\emptyset$.

Suppose now that $p\geqslant3$ and $n_1+..+n_{p-1}< n_p$. Put
$N_1=n_1+\cdots+n_{p-1}$. In this case, since the last
component $\fk_p\simeq\mathfrak{u}(n_p)$ of
$\fk$ is dominant in $\fk$, the calculation problem can be
reduced to the previous case with
$n_p=n_1+\cdots+n_{p-1}$. To this end we consider the
representation of the Lie group
$K_p\subset K$ with the Lie algebra
$\fk_p\simeq\mathfrak{u}(n_p)$ in the
$\fk$-submodule $V=\sum_{j=1}^{p-1}V^{j,p}$ of
$\fm$. Identifying $V$ with the space $M_{N_1,n_p}$ of complex
$N_1\times n_p$ matrices, the
$\Ad$-action of $K_p$ in
$V$ is described as follows: $k\cdot B=Bk^{-1}$, where
$k\in K_p=U(n_p)$, $B\in M_{N_1,n_p}$. Since the number of rows
$N_1$ in $B$ is less then its number of columns
$n_p$, then the $\Ad(K_p)$-orbit of $B$ in
$V$ contains a matrix in which the last
$n_p-N_1$ columns vanish. In other words, each element of
$\fm$ is
$\Ad(K)$-conjugated to some element of the first component
$\mathfrak{u}(2N_1)$ of the Lie algebra
$\mathfrak{u}(2N_1)\oplus\mathfrak{u}(n_p-N_1)\subset\mathfrak{u}(n)$
and, consequently,
$\dim\fk^x\geqslant 1+\dim\mathfrak{u}(n_p-N_1)$
($\dim\fz(\fg)=1$) for any
$x\in R(\fm)$. But we showed above that for some element
$x_0\in\tilde\fm$ we have an equality, i.e. this element
belongs to the set $R(\fm)\cap\tilde\fm$. Since the set
$R(\tilde\fm)$ is Zariski open in
$\tilde\fm$, we can suppose without losing generality that
$x_0\in R(\fm)\cap R(\tilde\fm)$ and
$\fk^{x_0}=\fz(\fg)\oplus\mathfrak{u}(n_p-N_1)$.
Taking into account that
$\mathfrak{u}(2N_1)\oplus\mathfrak{u}(n_p-N_1)$
is a maximal subalgebra of
$\mathfrak{u}(n)$ (the pair
$\bigl(\mathfrak{u}(n),\mathfrak{u}(2N_1)\oplus\mathfrak{u}(n_p-N_1)\bigr)$
is a symmetric pair), we obtain that the centralizer
$\fg_0$ of $\fk^{x_0}$ in $\fg$ is the algebra
$\mathfrak{u}(2N_1)\oplus\fz(\fg)$ containing the element
$a$. Now applying Theorem~\ref{th.13} and Remark~\ref{re.12}
we reduce our problem to the case when
$n_p=n_1+\cdots+n_{p-1}$. In other words,
$O^\mathrm{Kr}(\fm)\cap\tilde\fm\ne\emptyset$,
i.e. the condition of Theorem~\ref{th.10} holds.

Thus the new proof of Theorem 4 in~\cite{DGJ09} is obtained, i.e. the
following theorem is proved.
\begin{theorem}\label{th.16}~{\upshape\cite{DGJ09}}
We retain the notation of Theorem~\ref{th.10}. Let
$G=U(n)$, $\tilde G=SO(n)$ and
$\tilde K=SO(n_1)\times..\times SO(n_p)$,
$n=n_1+..+n_p$, with the standard embedding of
$\tilde K\subset\tilde G\subset G$. The set of functions
$\mathcal{F}(\fg,\tilde\fm)$ is a complete involutive subset
of the Poisson algebra
$(A^{\tilde G},\tilde\eta^\mathrm{can})$. The set of functions
$\mathcal{H}(\tilde\fg^*)\cup \mathcal{F}(\fg,\tilde\fm)$
is a complete involutive subset of the algebra
$(\mathcal{E}(T(\tilde G/\tilde K)),\tilde\eta^\mathrm{can})$.
The Hamiltonian flow with the Hamiltonian function
$\tilde H_{a,b}$ on $T(\tilde G/\tilde K)$
is completely integrable by means of real-analytic integrals
$\mathcal{H}(\tilde\fg^*)\cup \mathcal{F}(\fg,\tilde\fm)$.
\end{theorem}

\begin{remark}\label{re.17}
Let $\tilde\fg=\mathfrak{so}(n)$ and let
$\tilde\fk=\mathfrak{so}(n_1)\oplus\cdots\oplus\mathfrak{so}(n_p)$.
Suppose that $n_p\leqslant n_1+..+n_{p-1}$. Then the set
$\tilde\fm=\tilde\fk^\bot$
($\tilde\fg=\tilde\fm\oplus\tilde\fk$)
contains regular elements of the Lie algebra
$\mathfrak{so}(n)$ (see for example the proof of
Theorem~\ref{th.16} above). In~\cite{DGJ09x}
Dragovi\'c,  Gaji\'c and Jovanovi\'c
remarked that the first part of the proof of their Theorem 4
published in~\cite{DGJ09} is not complete. Namely,  relation (29)
of~\cite[Remark 1]{DGJ09}, i.e., the completeness of some family of integrals,
holds for all elements from
some open dense subset $O$ of $\mathfrak{so}(n)$
(for the so-called generic elements).
However, although
$\tilde{\mathfrak m}$ contains regular elements of
$\mathfrak{so}(n)$,  relation (29) holds on
$\tilde{\mathfrak m}$, i.e. $O\cap\tilde{\mathfrak m}\ne\emptyset$,
if and only if all $n_1,\dots,n_p$ are $\leqslant 2$
(see~\cite[p.~1287]{DGJ09x}). Therefore,
relation (36) of~\cite{DGJ09},
and, consequently, the first part of the proof of Theorem 4~\cite{DGJ09}
needs an additional argumentation in the case when  $n_p\geqslant3$.
\end{remark}

\end{document}